\documentclass[12pt]{amsart}
\usepackage{fullpage,amsmath,amsthm,amsfonts,amssymb,graphicx,amscd,float,hyperref,
setspace}
\usepackage[shortlabels]{enumitem}
\usepackage{color, xcolor}
\usepackage{tikz-cd}
\usepackage{xypic}
\usetikzlibrary{cd}
\usepackage{wrapfig}
\title{KarenAI}

\newtheorem{thm}{Theorem}[section]

\newtheorem{lem}[thm]{Lemma}

\newtheorem{prop}[thm]{Proposition}

\newtheorem{cor}[thm]{Corollary}

\theoremstyle{definition}

\newtheorem{df}[thm]{Definition}
\newtheorem{conjecture}[thm]{Conjecture}
\newtheorem{rmk}[thm]{Remark}

\newcommand{\mA}{\mathcal{A}}

\newcommand{\mG}{\mathcal{G}}

\newcommand{\mU}{\mathcal{U}}

\newcommand{\mW}{\mathcal{W}}

\newcommand{\LAU}{\mathfrak{u}}

\newcommand{\CC}{\mathbb C}

\newcommand{\RR}{\mathbb R}

\newcommand{\from}{\colon}

\newcommand{\im}{\text{Im}}

\newcommand{\gl}{\mathfrak{gl}}

\DeclareMathOperator{\Tr}{Tr}
\newcommand{\TrR}{\Tr_{\mathbb{R}}}
\theoremstyle{definition}

\theoremstyle{plain}

\newtheorem*{mainthmA}{Theorem A}
\newtheorem*{mainthmA'}{Theorem A$'$}
\newtheorem*{mainthmB}{Theorem B}
\newtheorem*{mainthmC}{Theorem C}
\newtheorem*{mainthmD}{Theorem D}
\newtheorem*{mainthmE}{Theorem E}
\newtheorem*{mainthmF}{Theorem F}
\newtheorem*{mainpropB'}{Proposition B$'$}

\newtheoremstyle{TheoremNum}
{8.0pt plus 2.0pt minus 4.0pt}{8.0pt plus 2.0pt minus 4.0pt} %%% space between body and thm
{\itshape} %%% Thm body font
{-0.15cm} %%% Indent amount (empty = no indent)
{\bfseries} %%% Thm head font
{.} %%% Punctuation after thm head
{ }  %%% Space after thm head
{\thmname{#1}\thmnote{ \bfseries #3}}%%% Thm head spec
\theoremstyle{TheoremNum}

\begin{document}

\title{The Continuum Limit of the Deep Linear Network}
\author{Catherine Eva Pfaff}
\author{Karen Uhlenbeck}

\address{{\tt Center for Communications Research - Princeton}}
\email{catherine.pfaff@gmail.com}

\address{{\tt Institute for Advanced Study}}
\email{uhlen@ias.edu}

\date{}
\maketitle

\begin{abstract}
The Deep Linear Network (DLN) is a simple version of a Multilayer Perceptron (MLP) in which the activation function is the identity and the bias vector is the zero vector. We examine the behavior of the model when the number of weight matrices, which is finite in the original model, goes to infinity. We show that a plausible continuum limit is a dimensionally reduced version of the Yang-Mills equation on the torus, revealing an unexpected geometric interpretation of infinite-depth DLNs. While we do not prove this rigorously, it is a very pretty calculation. It also gives some insight into several features of the construction. We note that the choices of boundary conditions on the ends of the discrete series survive  better in the continuum limit. Also, the regularizing functional which uses the moment map survives nicely in the continuum limit, whereas the regularizing functional using the size of approximating terms has no sensible limit. This surely is useful information for the choices we make in the finite-dimensional theory, and might also be useful in the more complicated MLP setting.
\end{abstract}

\section{Introduction}\label{s1:Intro}

In October 2025 Govind Menon gave a talk to mathematicians on the geometry of a Deep Linear Network (DLN). In MLP language, a DLN is an MLP where each activation function is the identity function and each bias vector is the zero vector. In practice, this means that a DLN operates on its input by applying a product matrix $X=W_N\cdots W_1$. To read this paper, one does not have to be familiar with the general model.

An interest in DLNs stems from the redundancy of the many choices of the $W_p$ leading to the same $X$. This redundancy leads to a geometric structure that is both interesting in its own right and can inform training dynamics. Saxe, McClelland, and Ganguli \cite{SMG14} showed that, despite the linearity of the networks themselves, DLNs have nonlinear gradient descent dynamics on weights that depend on the depth of the network. Further, they show that the training dynamics exhibit some of the phenomena observed in nonlinear neural networks. Perhaps even more interestingly for what we do, Arora, Cohen, and Hazan \cite{ACH18} were able to show that increasing depth in a DLN can accelerate gradient descent under certain assumptions.

We explore the continuum limit of this geometric structure here. In particular, Lindsey and Menon's crediting of the Kempf-Ness theorem as an essential tool for studying a DLN brings to mind the application to the heat equation for Yang-Mills on surfaces, which uses an infinite dimensional version of the theorem.  
Looking at the continuum limit of the DLN, the symmetries are identical in structure to the real and complex gauge changes for Yang-Mills on a Riemann surface. It was somewhat surprising that the theory turned out to be equivalent up to normalization (with a little fiddling of the gauge at the boundary) to a dimensionally reduced Yang-Mills theory from a torus to a circle with an additional reflection to fix the boundary conditions. We note that in \cite{cohen2023deep} Cohen, Menon, and Veraszto also take a continuum limit. However, our work is different from theirs in being focused on the continuum limit of the gauge-theoretic structure and its relationship to dimensionally reduced Yang–Mills theory.

What we indicate in this paper is that a dimensionally reduced Yang-Mills theory on the torus with a reflection is a plausible $N\to\infty$ limit of the DLN when $N$ is the number of layers in the network. The result obtains ODE's on the line from the limits of discrete versions. The discrete objects are obtained from continuous ones by identifying the $p/N$ on the line with the index-$p$ discrete object - the limits as $N\to\infty$ are then our continuous theory. We believe that a rigorous proof would be pretty straightforward, but lengthy, so we have not included it in this introductory paper.

In all theories which construct algorithms to identify information from new data, the aim in training is to determine the method for finding the best fit given large amounts of training data. Usually this involves a large number of intermediate steps. In the DLN we are starting from, this training flow occurs on $N$  $n\times n$ weight matrices $W_p \in GL(n,\CC)$. The product $X(1) =  W_N W_{N-1}\cdots W_p\cdots W_1$  represents the final information we want, and the $W_p$ are the intermediate steps. To the uninitiated, this seems quite strange, but rest assured it works. Because only $X(1)$ involves the final data, we have a geometry for the intermediate $W_p$. We may multiply a $W_p$ by an element of $GL(n,\CC)$ on the right as long as we multiply $W_{p-1}$ by the inverse of that element. This constructs a large ``gauge group'' which in this paper we always call the complex gauge group. When we only use unitary matrices, the group we use not only preserves the endpoints, $X(1)$ on the left and the identity on the right, but turns out to preserve all the geometry we construct.  We call this ``the gauge group''. It acts exactly like the gauge group in the physics and mathematical Yang-Mills models.  We know of no other theory that has such a gauge group, and find this very interesting.  The complex gauge group does appear in the study of Yang-Mills on a Riemann surface, which is where the idea for this paper came from.

The direction of the regularizing flow, which uses only the complex gauge group, is added to the more or less usual gradient flow which acts on the $X(1)$ term. This simplifies the weight matrices $W_p$ without affecting the flow on $X(1)$. The finite theory uses several regularizing flows, which we discuss in detail in Section \ref{s4:GeometryGrpOrbitsAndRegularizingFlows}.  It turns out that most of them (we don’t pretend to know all) flow to the zeros of the moment map for the gauge group actions. Hence the moment map is important in both the discrete and limiting theories. The moment map is reviewed in Section \ref{s2:LDNandContinuumLimit}  for those who are unfamiliar with symplectic geometry.

We consider two regularizing flows on the Hilbert space of connections. One is the standard Yang–Mills heat flow, arising as the gradient flow of the Yang–Mills functional. 
This turns out not to have viable boundary conditions when restricted to $\left[0,1\right]$.
The second is a somewhat unusual flow which can be interpreted as a smooth vector field in a Hilbert space of connections. 

The outline of the paper is as follows. Section \ref{s1:Intro} is the introduction. At the end of this section we give some comments on what one might learn about the finite theory from the infinite theory.  

Section \ref{s2:LDNandContinuumLimit} reviews step by step the geometry for finite $N$ and derives the continuum limit from the (usual) finite theory. We know the finite theory based on $N$ $n\times n$ invertible matrices $W_p$ acted on by $N - 1$ matrices $G_p$.

It is not necessary to explain here in detail how the discrete theory limits on a continuous theory, as that is carefully done in Section \ref{s2:LDNandContinuumLimit}.  It won’t hurt to give a few examples. Given the continuous (as smooth as we need) map $a\from \left[0,1\right]\to \mathfrak{gl}(n,\CC)$, we construct the discrete $W_p$ by setting $W_p = I + \frac{1}{N}a(\frac{p}{N})$. The gauge group is simple. 
Given $g\from \left[0,1\right]\to GL(n,\CC)$, we set $G_p = g(\frac{p}{N})$. We then find that, in the limit, the limiting theory from the finite theory has a connection $a$, with the standard action by a gauge group element g on it. We set $X_p=X(\frac{p}{N})=W_pW_{p-1}\cdots W_1$, which just yields in the limit a section of the trivial bundle that $g$ and the connection ($\frac{\mathrm{d}}{\mathrm{d}x}-a$) act on. The moment map, which is again described later,  for the discrete model converges to the moment map (curvature for Yang-Mills!) for the continuous model. Only after we construct the theory do we note which Hilbert spaces are natural in the problem. A comment here is that Yang-Mills theory is beautifully consistent up to dimension 4, so the analysis for the 1-dimensional problem is on the simple side.

The following two theorems (given in Section \ref{s2:LDNandContinuumLimit}) describe a plausible continuum model and geometry of the continuum limit. We use the notation $H^k$ for the Sobolev class
$$H^k := \{\text{functions } f \mid \text{each (weak) derivative } f^{(n)} \text{ with } 0\leq n\leq k \text{ is square integrable}\},$$
where the domain and range vary by circumstance.

Let $W:=\left(W_1,\dots, W_N\right)$.\\

\begin{mainthmA}[Theorem \ref{prop2.1}]\label{mainthmA}
We choose $a\from \left[0,1\right]\to \mathfrak{gl}(n,\CC)$, and $X\from[0,1] \to GL(n,\CC)$, and a gauge element $g$, and find $W_p$, and $G_p$, and $X_p$, as above. 
As $N\to\infty$, (a)-(d) below are plausible limits. 
\begin{itemize}
  \item[(a)] The continuum limit $X\from[0,1] \to GL(n,\CC)$ satisfies $\frac{\mathrm{d}}{\mathrm{d}x} X = a X$, $X(0) = I$.\\
      
  \item[(b)] The limiting gauge group $\mG^{\infty} ~:=~\{g\from \left[0,1\right]\to GL(n,\CC)\mid g\in H^2 \text{, } g(0) = g(1) = I\}$ acts on the connection $a$ from part (a) by $$(g, a)\mapsto \left(\frac{\mathrm{d}}{\mathrm{d}x} g + g a\right)g^{-1}.$$
      
  \item[(c)]  The tangent space at a to the orbit of the limiting gauge group $\mG^{\infty}$ is the image of
       $${ \frac{\mathrm{d}}{\mathrm{d}x} \alpha - \left[a,\alpha\right]\quad \text{for}\quad \alpha(0) = \alpha(1) = 0}.$$         
  \item[(d)] The moment map $\omega$ for the action of the unitary gauge group \\ $\mU^{\infty} ~:=~\{g\in\mG^{\infty}\mid g(x)\in U(n)\text{ for all } x \}$ is  
             $$\omega ~ = ~ \frac{\mathrm{d}}{\mathrm{d}x}\left(a + a^\ast\right) - \left[a,a^\ast\right].$$  \\
\end{itemize}
\end{mainthmA}

\begin{mainthmA'}[Theorem \ref{prop2.1'}]\label{mainthmA'}
The limiting configuration space 
$$\mA = \{ a \from [0,1]\to \mathfrak{gl}(n,\CC) \mid a\in H^1\}$$ 
carries the following structures:\\
\begin{itemize}
  \item[(a)] $\mA$ is a linear space, with 
  
       $$\langle a,b\rangle ~ = ~ \int_0^1  \TrR a(x)b(x)^\ast \mathrm{d}x.$$ \\
The symplectic structure is $\langle ia, b\rangle$. \\
      
  \item[(b)] The gauge group $\mU^{\infty}$ of unitary matrices in $\mG^{\infty}$ leaves the geometry invariant, whereas $\mG^{\infty}$ leaves only $X(0) = I$ and $X(1)$ invariant.\\
         
  \item[(c)] The regularizing function $||W||^2$ becomes infinite in the limit. \\
\end{itemize}
\end{mainthmA'}

In Section \ref{s3:DimReducedYangMillsFlow}, we describe the relationship between our continuum limit and Yang-Mills on the torus. The Yang-Mills theory on the torus $\left(x,y\right)$ uses the product connection $\left(\left(\frac{\mathrm{d}}{\mathrm{d}x} - u\right), \left(\frac{\mathrm{d}}{\mathrm{d}y} + ih\right)\right)$, with $u$ and $ih$  skew-Hermitian. Dimension reduction in the $y$ variable means that $u$ and $ih$ do not depend on $y$. It does not mean discarding the $ih$ variable, which is often called the Higgs field in the similar higher dimensional examples in geometry and physics. In our theory, we write the complex connection $a$ as $a = h + u$ (the reason for our weird inclusion of the $i$ in the original description), where $h = \frac{1}{2}\left(a + a^{\ast}\right)$ and $u = \frac{1}{2}\left(a - a^{\ast}\right)$. To reduce the problem on the circle to a line, we add the reflection $\left(a(x)\right)^{\ast}=a(-x)$ and $\left(g(x)\right)^{\ast}=\left(g(-x)\right)^{-1}$. Here we run into problems, as $u$ behaves badly and the unitary gauge groups are different.  We describe how both problems are due to differences in the (unitary) gauge group.  
In the reflected model, when $u$ has this jump, there is no difficulty with the flow, but the jump is preserved. In the end this might contain information useful in constructing the numerics for the flows.

The main result of Section \ref{s3:DimReducedYangMillsFlow} is Theorem \ref{thm3.3}. We restate the gist of it here as Theorem B.

\begin{mainthmB}[Theorem \ref{thm3.3} abbreviated]\label{mainthmB}
With the above identifications, Yang–Mills theory on the torus reduced by setting terms independent of $y$ to the circle, and further reduced to [0,1] by the above reflections, is equal to the theory for our infinity limit, with slight modifications in the symmetry groups.
\end{mainthmB}

However, when we look at the heat equation for Yang-Mills, we find it difficult to use on the interval $\left[0,1\right]$ without adding extra boundary conditions. This gives us Proposition B$'$, which is Proposition \ref{propB'}.

\begin{mainpropB'} The usual gradient descent for Yang-Mills is a heat equation, which does not restrict to the interval $\left[0,1\right]$ to give a gradient descent on $\left[0,1\right]$. 
\end{mainpropB'}

Section \ref{s4:GeometryGrpOrbitsAndRegularizingFlows} treats the gradient flows, first for finite $N$, and then for the limit as $N\to\infty$. We repeat the calculations for finite $N$ first as these calculations are more complicated. Two of the gradient flows for the two main regularizing functionals we know about turn out to be the same! 
In the continuum limit, only one gradient flow survives.

We write the heat equation for Yang-Mills in terms of the dimensionally reduced variables:

\begin{mainthmC}[Proposition \ref{prop3.2} and following remarks]\label{mainthmC}
In terms of dimensionally reduced variables
         $$\frac{\mathrm{d}}{\mathrm{d}s} h ~ = ~ 2~  \left[~\frac{\mathrm{d}}{\mathrm{d}x} - u,~ \left[\frac{\mathrm{d}}{\mathrm{d}x} - u ,~ h\right]~\right],$$
$$\frac{\mathrm{d}}{\mathrm{d}s} u ~ = ~  -2 ~ \left[~h,~ \left[ \frac{\mathrm{d}}{\mathrm{d}x} - u,~h\right]~\right].$$
\end{mainthmC}

\begin{mainthmD}[Theorem \ref{prop5.6}, Corollary \ref{cor5.7}]\label{mainthmC}
The rescaled limit as $N \to \infty$ of our first gradient flow, which we will call the group geometric flow, is the dimensionally reduced heat equation from Yang-Mills.

Hence, by Proposition B$'$, it is not converging to a flow with viable boundary conditions.
\end{mainthmD}

The continuum limit of the second flow was harder to find. Our first attempt turned out to be incorrect. When we did find it, it was quite interesting.

\begin{mainthmE}[Theorem \ref{thm5.12}, Corollary \ref{corofthm5.12}]\label{mainthmE}
The continuum limit of the linear geometric flow
$$\frac{1}{2}\frac{\mathrm{d}}{\mathrm{d}s}a ~=~ v-a,$$
where
$$v-a ~=~ \left[~\left(\frac{\mathrm{d}}{\mathrm{d}x}-a\right),~\tau~\right],\quad\quad
\tau(0)=\tau(1)=0$$
$$\left[~\frac{\mathrm{d}}{\mathrm{d}x}+a^{\ast},~ \left(v+a^{\ast}\right)~\right] ~=~ 0.$$
Here $v\in H^1$ is smoothly dependent on $a\in H^1$, and $\frac{\mathrm{d}}{\mathrm{d}s}a$ is a smooth vector field on $a\in H^1$.
\end{mainthmE}

We conjecture (Conjecture \ref{conj1}) that the projected flows in the finite-dimensional model converge, as $N \to\infty$, to the corresponding projected flows in the continuum model.

Section \ref{s5} completes the proof that the conjectured limits of the loss function we aim to minimize and of the regularizing functional behave as we would like them to behave. This is Theorem \ref{thm5.5new}. We reformulate this result in less technical language here.

\begin{mainthmF}[Theorem \ref{thm5.5new} abbreviated] 
The downward gradient flow for $F(a)+\nu R(a)$, 
where $F(a) = f(X(a))$ and $R(a) = \langle\omega,~\omega\rangle$, is
$$
\frac{\mathrm{d}a}{\mathrm{d}s} ~=~ 
-\left(X(x)^{-1}\right)^{\ast}X(1)^{\ast}
\left(\nabla_Xf\right)X(x)^{\ast}+
\nu\frac{\mathrm{d}a}{\mathrm{d}s},
$$
where $\frac{\mathrm{d}a}{\mathrm{d}s}$ is a downward gradient for $\langle\omega,~\omega\rangle$ leaving $X(1)$ fixed. Then 
$$\frac{\mathrm{d}}{\mathrm{d}s}F(a) ~=~ 
-dF(a)\left(\nabla F \right)
 ~=~ 
-\left\langle \nabla F, \nabla F \right\rangle,$$
$$\frac{\mathrm{d}}{\mathrm{d}s}R(a) ~=~ 
2\nu \left\langle \omega, \frac{\mathrm{d}}{\mathrm{d}s}\omega \right\rangle.
$$\\
For the flow described in Theorem E, $\frac{\mathrm{d}}{\mathrm{d}s}R(a)=-8\nu \left\langle \omega, \omega \right\rangle$.
\end{mainthmF}

We are new to this field, so have little to say about where it is going. The first point is that the norm squared of the moment map survives as a good regularizing functional in the limit. The norm squared of the $W_p$, or any other norm, will only become infinite. Even in the finite case, it seems to make sense numerically to deal with the smaller numbers coming from $\Omega$ rather than the larger numbers coming from $X$. 

One of the main technical points is in handling the boundary term corresponding to $W_1$ and $W_N$.  To use the gradient flow, it is necessary to use a (unitary) gauge symmetry to regularize the boundary conditions. It might be helpful to take this into account when the gradient flows are constructed numerically. 

A further point is that any choice of gradient flow for the moment map squared is a good choice for a regularizing flow. This does not seem to be the case for the flows coming from norms of the $W=\left(W_1,\dots, W_N\right)$, unless they by accident happen to be gradient flows for the moment map. Hence, we believe the moment map is better to work with than the $W$'s. Since the gradient flows we use for the norm squared of the $W$'s are also gradient flows for the norm squared of the moment map, this should be automatic. 

We found two regularizing flows associated to the moment map. The first was a heat equation, which is ill posed with the desirable boundary condition on $\left[0,1\right]$. The second one converges. The flows labeled the linear geometric flows in Definition \ref{linear_geometric_flow} are better suited for large $N$ than those labeled the group geometry flows in Definition \ref{group_geometry_flow}. Finally, the description given for points at which the moment map vanishes  does not at all carry over to the continuum limit. This should be considered in describing this submanifold. 

The relationship with Yang-Mills is reminiscent of the appearance of $KdV$ in many different areas. Material can be found in the papers of Menon, both solely authored and with collaborators. The standard references for gauge theory that we know of are \cite{baez1994gauge, donaldson1997geometry}. A more extensive treatment of related aspects of the infinite-depth geometry of deep linear networks appears in the doctoral thesis of Tejas Kotwal \cite{kotwal-thesis}.

An interesting direction we did not pursue here is understanding how including nonlinear activation functions (and even bias vectors) changes the underlying geometry and training dynamics of the network.

\subsection*{Acknowledgements}

The authors are grateful to Govind Menon for introducing the problem and the hospitality of the Institute for Advanced Study, without which this paper would not have come about. The first author is grateful for the opportunity to play an apprentice role in this work. The second author is grateful to Linda Ness for her interest in their work. We thank Tejas Kotwal and Laura Fredrickson for their comments and improvements to our paper. We also acknowledge use of ChatGPT in the revision process.

\section{The Linear Deep Learning Network and the continuum limit}
	\label{s2:LDNandContinuumLimit}  

For the purposes of training a neural network, the goal of the subject is to find a good description of the minima of a function $f\from GL(n,\RR) \to \RR$, that will be a loss function for training this neural network. The DLN at level $N$ (i.e. that with $N$ layers) does this by writing  $X = X(W) =  W_N \cdots W_p \cdots W_1$ where each weight matrix $W_p\in GL(n,\RR)$, and $W = \left(W_1, \dots,W_N\right)$, and $F(W) = f(X(W))$. However, the geometry is based on $W_p\in GL(n,\CC)$, as is the connection with Yang-Mills, and we develop the geometry from that point of view. In applications, we would simply apply our results to the real matrices, but to understand the geometry we must use all of $GL(n,\CC)$. We study the gradient flow of $F$, but we add to it a regularizing term $\nu R(W)$ which picks out of the possibilities of describing $X$ those which are ``simple''.  We treat these flows separately in Section \ref{s4:GeometryGrpOrbitsAndRegularizingFlows}.

We highlight that our definition of the DLN at level $N$ is based on $N$ elements of the Lie group $GL(n,\CC)$ rather than $GL(n,\RR)$, and we denote those also by $W_p$. We view $GL(n,\CC)$ as an open subset of the ambient linear space $\mathfrak{gl}(n,\CC)$, and use the linear geometry it inherits as such. We denote the parameter space with parameters $W_1, \dots, W_p, \dots, W_N$ by $\mW$. 
Note the typical redundancy in the number of parameters. Each $W_p$ contributes $n^2$ weights to be trained, so that there are $N\times n^2$ weights in total. We use those to describe $n^2$ values in $X(1)$.

One looks for the minimum of the loss functional, which depends on $X = W_N\cdots W_1$ for $W = (W_1, \dots,W_N)\in\mW$. This is done by standard gradient descent. However, as explained in the introduction, there is a group of gauge transformations $\mG$ acting on $\mW$, which preserve  $X(W)$, and hence $f(X(W))$. The minimization procedure is to add to the function  $f(X(W))=F(W)$ the regularizing term $\nu R(W)$, whose gradient flow lines are infinitesimal gauge transformations, and hence the flows do not change the value of $X(W)$.  In other words, the flow in the quotient space $\mW$/$\mG$ does not change, but the representation by $W$ does change.  We are not interested in this paper in the function $F$, but in the structure of the action of the gauge group $\mG$ on $\mW$ and the possible regularizing terms. We only briefly describe the structure for finite $N$, as the details are easily found in other references.

The complex gauge group is 
$$\mG = GL(n,\CC)\times\cdots\times GL(n,\CC)  = \{ \left(G_1, \dots,  G_p, \dots, G_{N-1}\right)~:~ G_p\in GL(n,\CC)\}$$ 
and the real (unitary) gauge group is the subgroup $\mU = U(n)\times\cdots\times U(n)$ of $\mG$ for which each $G_p$ is unitary. The real gauge group acts as symmetries of the entire geometry, whereas the complex gauge group changes the geometry, while still preserving $X(W)$.  In the Yang-Mills setting, the real gauge group consists of what we usually think of as gauge transformations (which yields an infinite dimensional quotient space). The complex gauge group appears only for Riemann surfaces and preserves only the complex structure (which yields a finite dimensional quotient space).

\begin{df}\label{df2.1}
The action $\mG \times \mW \to\mW$ is given by
$$
\begin{cases}
W_1\mapsto G_1 W_1, \\
W_p \mapsto G_p W_p (G_{p-1})^{-1}\quad \text{ for } p\neq 1,N,\\
W_N \mapsto W_N(G_{N-1})^{-1}.\\
\end{cases}
$$
\end{df}

In algebraic geometry, the quotient space $\mW$/$\mG$  is not necessarily Hausdorff. Unstable orbits have to be removed to define the quotient, although there are no unstable orbits in our circumstance. This is where the Kempf-Ness theory applies, the squared norm of the moment map is the classical Kempf–Ness functional (see, for example, Kempf–Ness \cite{kempfness1978} or Kirwan \cite{kirwan2020}). There is a general theory that is applicable to the complex setting and gives information on the restriction to the real case.

Since $\mW$ is an open subset of a complex vector space, it carries the standard symplectic structure induced by its complex structure. Let $\LAU$ denote the Lie algebra of $\mU$ and let $\LAU^\ast$ be its dual. Recall that the moment map for a Lie group action is a map $m\from \mW \to \LAU^\ast$ such that for each $A\in\LAU$, the corresponding Hamiltonian function $m(A)$ on $\mW$, equipped with its symplectic structure, generates the infinitesimal action of the Lie algebra $\LAU$ on $\mW$. 

Before describing the moment map in our setting, we provide some background, for those unfamiliar with the subject. The symplectic gradient is the complex number $i$ times the real gradient if the symplectic structure comes from a complex structure. Also, one often identifies the dual of the Lie algebra with the Lie algebra by using the inner product  $\langle A,B\rangle ~=~\TrR\left(AB^\ast\right)$, where $\TrR$ denotes the real trace. However, in the literature we follow, we identify the dual of the Lie algebra with the symmetric matrices by using the symplectic inner product $\im \left(\Tr \left(AB^\ast\right)\right)~=~ \langle A,iB\rangle$.  

In our setting, independent of the ambiguities explained in the previous paragraph, the moment map for this problem is given by
$$m(W) = \left(\Omega_1, \dots, \Omega_p, \dots, \Omega_{N-1}\right),$$
where $\Omega_p = W_{p+1}^\ast W_{p+1} - W_p W_p^\ast$ for each $p\in\{1,\dots N-1\}$. Note that $\Omega_p = \Omega_p^\ast$ and let $\Omega$ denote the moment map $\left(\Omega_1, \dots, \Omega_p, \dots, \Omega_{N-1}\right)$. It may be noteworthy that setting $\Omega=0$ is a precondition of \cite{ACH18}. This is different from adding a regularizing flow to training, which is what we do.

We repeat some of what is in the literature on the problem for finite $N$. Two relevant surveys are \cite{m23} and \cite{wendin-altafini}.  The Riemannian geometry and gradient flows associated with DLNs are also studied in \cite{brtw22}.

\begin{prop}[see \cite{m23}, \cite{km25}]\label{prop2.2}
The following are properties for finite N:\\
\begin{itemize}
  \item[(a)]	$\mW$ is an open subset of a (complex) linear space with an inner product 
        $$\left\langle W,V\right\rangle ~ = ~ \sum_p \TrR W_p V_p^{\ast}.$$  
        The symplectic structure is $\left\langle iW,V\right\rangle$.\\
  \item[(b)] $X\left(p/N\right) = W_p X\left(\left(p-1\right)/N\right)$ where $X\left(p/N\right) = W_p\cdots W_1$.\\
  \item[(c)] The gauge group $\mG$ acts on $\mW$ by
  $$
\begin{cases}
W_1\mapsto G_1 W_1, \\
W_p \mapsto G_p W_p \left(G_{p-1}\right)^{-1}\quad \text{ for } p\neq 1,N,\\
W_N \mapsto W_N\left(G_{N-1}\right)^{-1}.\\
\end{cases}
$$
  \item[(d)] The gauge group of unitary matrices $\mU$ leaves the geometry invariant, whereas $\mG$ leaves only $X = X(1)$ (and $X(0) = I$) invariant.\\
  \item[(e)] The moment map for the action of $\mU$ on $\mW$ is $\Omega$.\\
  \item[(f)] The tangent space to the orbit of $\mG$ on $\mW$ is               
 $$\delta W_p~=~ \delta G_p G_p^{-1} W_p - W_p \delta G_{p-1}G_{p-1}^{-1}.$$
  \item[(g)] Several regularizing functions  $\left\langle W,~W\right\rangle$ and $\left\langle \Omega,~\Omega\right\rangle$ are used.
\end{itemize}
 \end{prop}

We first explain the notation $\delta G_p$, which is shorthand for an element in the tangent space at $G_p$: $\delta G_p = \frac{\mathrm{d}}{\mathrm{d}s} G_p$.  Since we never use the parameter $s$, we omit it.  If $y = \LAU(x)$, we write $\delta y  =  \frac{\mathrm{d}}{\mathrm{d}x} \LAU \delta x$.  This is a classical notation which is not used much today. However, the more abstract definitions are wordy, and the basic ideas become obscure. 

Finally, note that $\alpha = \delta G_p G_p^{-1}$ is the Lie algebra element associated to the variation $\delta G_p$ of $G_p$.

\begin{df}\label{df2.3}
Start with $a\from \left[0,1\right]\to \mathfrak{gl}\left(n,\CC\right)$ and $g\from \left[0,1\right]\to GL\left(n,\CC\right)$.
Construct a finite model with $N$ terms by setting

\begin{itemize}
  \item[i.] $W_p = I + 1/N a\left(p/N\right)$ and 
  
  \item[ii.] $G_p = g\left(p/N\right)$.
\end{itemize}

We find it useful to let

\begin{itemize}  
  \item[iii.] $X\left(p/N\right) = W_p \cdots W_1$.
\end{itemize}
\end{df} 

To derive the model, we allow $a$ and $g$ to be as smooth as we need. We assume the notation of Definition \ref{df2.3}.

\begin{thm}[Continuum model]\label{prop2.1}
As $N\to\infty$, (a)-(d) below are plausible limits. \\
\begin{itemize}
\item[(a)] The continuum limit $X\from\left[0,1\right] \to GL\left(n,\CC\right)$ and $a\from\left[0,1\right] \to \mathfrak{gl}\left(n,\CC\right)$ satisfies $\frac{\mathrm{d}}{\mathrm{d}x} X = a X$, $X(0) = I$. Note that $X(0)$ is, as always, the identity or the theory will be trivial.\\
      
  \item[(b)] The gauge group $\mG^{\infty}$ from Definition \ref{df2.3}ii acts on the connection $a$ from Definition \ref{df2.3}i by 
      $$(g, a)\mapsto \left(\frac{\mathrm{d}}{\mathrm{d}x} g + g a\right)g^{-1}.$$
      
  \item[(c)] The tangent space at $a$ to the $\mG^{\infty}$-orbit consists of all vectors  
       $${ \frac{\mathrm{d}}{\mathrm{d}x} \alpha - \left[a,\alpha\right]\quad \text{with}\quad ~\alpha(0) = \alpha(1) = 0}.$$  
       
  \item[(d)] The moment map $\omega$ for the action of the unitary gauge group  $\mU^{\infty}$ is  
             $$\omega ~ = ~ m(a) ~ = ~ \frac{\mathrm{d}}{\mathrm{d}x}\left(a + a^\ast\right) - \left[a,a^\ast\right].$$    
      
      This is the limit of the discrete moment maps constructed from the formulas of Definition \ref{df2.3}.   \\
       
\end{itemize}
\end{thm}

In Corollary \ref{cor_for_2.4} we will show that the formulas we get for $a$, $g$, and $X$ are all consistent with $a\in H^1$, $g\in H^2$, $X\in H^2$, and $\omega$ in $H^0$. \\

\begin{proof}[Proof of Theorem \ref{prop2.1}]
We use smooth maps $a$, $X$, and $g$ to insure the above are well defined. 

Part (a) is obtained by taking the continuum limit of the recurrence for the partial products:

$$X\left(\frac{p}{N}\right) ~ = ~ W_p ~X\left(\frac{p-1}{N}\right) ~ = ~ \left(I + \frac{1}{N} a\left(\frac{p}{N}\right)\right) X\left(\frac{p-1}{N}\right).$$\\
Multiply by $N$ and rewrite this as

$$\frac{X_{p/N} - X_{(p-1)/N}}{1 ~/~ N} ~ = ~ a_{p/N}~ X_{(p-1) /N}.$$\\
In the limit we get (a). 

Similarly for (b): For $1<p<N$, start from the finite-$N$ gauge action to get

            $$W_p \mapsto  G_p W_p G_{p-1}^{-1}.$$
Thus              
          $$I + \frac{1}{N} a_{p/N} \mapsto$$ 
          $$g_{p/N}
          \left(I + \frac{1}{N} a_{p/N}\right) g_{(p-1)/N}^{-1} ~ = ~ 
          I+\frac{1}{N}\left(\frac{g_{p/N} - g_{(p-1) /N}}{1~/~N} + g_{p/N}~a_{p/N}\right)g_{(p-1)/N}^{-1}.$$ \\
 
Multiplying by $N$, replacing the difference quotient by a derivative, and then using the fact that $g$ is differentiable yields (b). 

Part (c) follows by differentiating the gauge action in (b) in the case that  $\delta g g^{-1} = \alpha$.

We are down to (d). First we show that the formula can be derived as a limit of the formula for large $N$.

      $$\Omega_p ~ = ~ W_{p+1}^\ast W_{p+1} - W_p W_p^\ast  ~ = ~ 
 \left(I + \frac{1}{N} a_{(p+1) /N}\right)^\ast \left(I+\frac{1}{N} a_{(p+1) /N}\right)  - 
           \left(I +\frac{1}{N} a_{p/N}\right)
           \left(I + \frac{1}{N}a_{p/N}\right)^\ast$$
$$  ~ = ~ \left(\frac{1}{N}\right)^2 \left(\frac{a_{(p+1) /N} + a_{(p+1) /N}^\ast - a_{p/N} - a_{p/N}^\ast}{1/N}   +
              a_{(p+1) /N}^\ast ~a_{(p+1) /N} - a_{p/N}~ a_{p/N}^\ast\right).$$

If we multiply by $N^2$, we get equation (d) in the limit. Thus $\Omega_p$ is of order $N^{-2}$. The unscaled moment map $\Omega_p$ therefore tends to zero in the continuum limit, while the renormalized quantity $N^2 \Omega_p$ converges to the continuum moment map $\omega$. This normalization by $N^2$ is the normalization underlying the continuum regularizer $||\omega||^2$.

Note that we are somewhat cavalier about keeping track of constants, which just determine the velocities of the gradient flows.  Also, we are not sure of the geometry behind why we need to multiply the limit in part (d) by $N^2$.  It must be a reflection of something going on in the geometry, just as $\langle W,~W\rangle$ in the limit becomes infinite and unusable.

Since the identification with Yang-Mills is not precise, we show here that $\omega$ is the moment map for the unitary gauge group action on $a$. Recall that we are using the symplectic structure to identify the dual of the Lie algebra. So we want the symplectic gradient of $\langle i\omega, \alpha\rangle$ to be the tangent vector to the orbit in the direction $\alpha\in H^2_0$. 

At $a$, this is just $\left[\left(\frac{\mathrm{d}}{\mathrm{d}x}-a\right),\alpha\right]$.
$$
\delta \left\langle i\omega, \alpha\right\rangle ~ = ~ 
\delta \left\langle i\frac{\mathrm{d}}{\mathrm{d}x}\left(a + a^{\ast}\right)
-\left[ a,a^{\ast}\right], 
\alpha\right\rangle ~ = ~ 
\left\langle i\left[\left(\frac{\mathrm{d}}{\mathrm{d}x}+ a^{\ast}\right),\delta a\right], \alpha\right\rangle+
\left\langle i\left[\left(\frac{\mathrm{d}}{\mathrm{d}x}-a\right),\left(\delta a\right)^{\ast}\right], \alpha\right\rangle.
$$
Because $\alpha ~ = ~ -\alpha^{\ast}$ and $\left[i\left(\frac{\mathrm{d}}{\mathrm{d}x}-a\right),\left(\delta a\right)^{\ast}\right]^{\ast} ~ = ~ 
-i\left[\left(\frac{\mathrm{d}}{\mathrm{d}x}+ a^{\ast}\right),\delta a\right]$, these terms are actually the same.

We now have
$$
2\left\langle i\left[\left(\frac{\mathrm{d}}{\mathrm{d}x}+ a^{\ast}\right),\delta a\right], \alpha\right\rangle ~ = ~ 
-2\left\langle i\delta a, \left[\left(\frac{\mathrm{d}}{\mathrm{d}x}-a\right), \alpha\right] \right\rangle.
$$

Here we are allowed to integrate by parts with no boundary term because $\alpha(0) = \alpha(1) = 0$. Now it is clear that the symplectic gradient is $-2\left[\left(\frac{\mathrm{d}}{\mathrm{d}x}- a\right),\alpha\right]$.

Note that the normalizations vary. We perhaps should have set $-\frac{1}{2}\omega$ as the moment map. The computation is more transparent if we identify the dual of the Lie algebra with the Lie algebra via the inner product, which makes $i\omega$ the moment map.
\end{proof}

\begin{cor}\label{cor_for_2.4}
All formulas of Theorem \ref{prop2.1} are consistent with 
$$g\in H^2, \quad X\in H^2, \quad a\in H^1, \quad \omega\in H^0.$$ 

In particular, 
      $\omega\in H^0$ when $a\in H^1$, and (d) turns out to work for $\alpha\in H^2$.
\end{cor}

\begin{proof}
(a) The analysis will work perfectly, as we can rewrite $\frac{\mathrm{d}}{\mathrm{d}x} X = a X$ as $\frac{\mathrm{d}}{\mathrm{d}x} X X^{-1} = a$.  Now $X$ and $X^{-1}$ are in $H^2$, $\frac{\mathrm{d}}{\mathrm{d}x} X$ is in $H^1$, and $H^2\cdot H^1\subset H^1$ as long as we are in dimensions less than 4! 

(b) Again, the Sobolev regularity is preserved.  We want $\left(\frac{\mathrm{d}}{\mathrm{d}x} g + g a\right)g^{-1}$ to be in $H^1$ when $a\in H^1$ and $g\in H^2$. This follows from several applications of the fact that $H^2\cdot H^1\subset H^1$ in dimension 1.

(c) We just need to know that the derivative of $a+a^\ast$ and the commutator $\left[a, a^\ast\right]$ are in $H^0$ when $a$ is in $H^1$.
Again, the analysis works perfectly.  

(d) If $g\in H^2$, then $\alpha=\delta g g^{-1}\in H^2$. Hence, $\frac{\mathrm{d}}{\mathrm{d}x} \alpha - \left[a,\alpha\right]\in H^1$ since $\frac{\mathrm{d}}{\mathrm{d}x} \alpha\in H^1$ and $\left[a,\alpha\right]\in H^1$.
\end{proof}

\begin{thm}[Geometry of the continuum limit]\label{prop2.1'}
The limiting configuration space $\mA$ carries the following structures:
\begin{itemize}
  \item[(a)] $\mA$ is a linear space, with 
         $$\left\langle a,b\right\rangle ~ = ~ \int_0^1  \TrR a(x)b(x)^\ast \mathrm{d}x.$$ \\
The symplectic structure is $\left\langle ia, b\right\rangle$. \\
      
  \item[(b)] The gauge group $\mU^{\infty}$ of unitary matrices in $\mG^{\infty}$ leaves the geometry invariant, whereas $\mG^{\infty}$ leaves only $X(0) = I$ and $X(1)$ invariant.\\
         
  \item[(c)] The regularizing function $||W||^2$ becomes infinite.  
Only $|| \omega ||^2$ survives as a regularizing function.\\
\end{itemize}
\end{thm}

\begin{proof}
Item (a) is mostly just a definition. Note, however, that this is an $H^0$ inner product and does not give the correct topology. This is the origin of the heat equation we will see in later sections. However, it does give the symplectic structure we use later. Note that in the approximation, for  $W = \left(W_1, \dots, W_N\right)$  with $W_p = I + \frac{1}{N} a\left(p/N\right)$
and $Y = \left(Y_1, \dots, Y_N\right)$ with  $Y_p = I + \frac{1}{N} b\left(p/N\right)$, we have that
$\left\langle W,~Y\right\rangle ~ = ~ \sum_1^N \TrR W_p Y_p^\ast$  is not related to the limit $\left\langle a, b\right\rangle$ in any meaningful way.  

We don’t really have to prove all of (b) as in the end the proof is contained in the literature on Yang-Mills, which we explain in Section \ref{s3:DimReducedYangMillsFlow}. We will see as we proceed that $\mU^{\infty}$ leaves the geometry invariant. To check the $\mG^{\infty}$ leaves $X(0) = I$ and $X(1)$ invariant, we need to check how $\mG^ {\infty}$ acts on $X$. 

We start with the answer.  We conjecture that $$g(X) = g X.$$ If this were true, then  
       $$\frac{\mathrm{d}}{\mathrm{d}x} \left(g X\right) \left(gX\right)^{-1} ~ = ~ \left(\frac{\mathrm{d}}{\mathrm{d}x} g\right)g^{-1} + g\left(\frac{\mathrm{d}}{\mathrm{d}x}X\right)\left( X^{-1} g^{-1}\right) ~ = ~  
\left(\frac{\mathrm{d}}{\mathrm{d}x} g + g a\right)g^{-1},$$ 
which gives  the correct multiplication. We could have started with this definition! Note that $H^2\times H^2$ lands in $H^2$ in dimensions less than 4, so we are on track with the analysis. 

It is clear that the finite-dimensional functional $||W||^2$ diverges in the continuum limit.  We are left with $||\omega||^2$ as the regularizing function.  Again, when we identify the $N = \infty$ limit with Yang-Mills,  $\omega$ is proportional to the curvature (under the standard identification), and then $\left\langle\omega, \omega\right\rangle$ is a scalar multiple of the usual Yang-Mills functional.

\end{proof}

\section{The dimensionally reduced Yang-Mills flow}
	\label{s3:DimReducedYangMillsFlow}

Before we move on to the flows, we explain the relationship with Yang-Mills. This will save us a certain amount of work, since the formalism is already set up. We first rewrite the above equations in terms of the Hermitian and skew-Hermitian parts of $a$.  

\begin{thm}\label{thm3.2}
Let  $a = h + u$ and $a^\ast = h-u$, so that $h=\frac{1}{2}\left(a+a^\ast\right)$ and $u=\frac{1}{2}\left(a-a^\ast\right)$. Then each of the following holds:\\
\begin{itemize}
\item[(a)] $\frac{\mathrm{d}X}{\mathrm{d}x} = \left(h + u\right)X$.\\
\item[(b)] 
$(g,h,u)\mapsto \left( \frac{\mathrm{d}}{\mathrm{d}x} g +g\left(h+u\right)\right)g^{-1}$ for each $g\in\mG^{\infty}$.
In particular, for each  $g\in\mU^{\infty}$, 
    $$u\mapsto \left( \frac{\mathrm{d}}{\mathrm{d}x} g +gu\right)g^{-1}
    \quad \text{and}\quad h\mapsto ghg^{-1}.$$
\item[(c)]	$\omega = 2\left[\left(\frac{\mathrm{d}}{\mathrm{d}x} - u\right),h\right]$.\\
\item[(d)]	If $\alpha = \alpha^\ast$ (so that $\alpha$ lies in the orthogonal space to the unitary Lie algebra), then the tangent vector generated by the infinitesimal gauge transformation $\alpha$ is
        $$\delta h ~ = ~ \left[\left(\frac{\mathrm{d}}{\mathrm{d}x} - u\right),~\alpha\right],$$
        $$\delta u ~ = ~ \left[\alpha,~h\right].$$
\end{itemize}
\end{thm}

\begin{proof}
We already explained the notation $\delta  v = \frac{\mathrm{d}}{\mathrm{d}s} v$.  This represents the tangent vector to the space that $v$ sits in at $v$, and behaves exactly as $\frac{\mathrm{d}}{\mathrm{d}s} v$.  Since we never actually use the parameter $s$, it is less confusing to avoid writing $s$. Note that we can also write $\alpha = (\delta g) g^{-1}$ as an element of the Lie algebra based at $g$.

Equation (a) is directly copied from Theorem \ref{prop2.1}(a). Equation (b) follows from the gauge action on $a=h+u$ and then, in the unitary case, decomposing the transformed connection into its Hermitian and skew-Hermitian parts.  
$$
\left(g,h,u\right) ~ = ~ 
\left(\left(\frac{\mathrm{d}}{\mathrm{d}x}g\right) g^{-1}+gug^{-1}\right) +
ghg^{-1}
$$
If $g^{\ast} ~ = ~ \left(g^{-1}\right)$, we get
$$
\left(\left(\frac{\mathrm{d}}{\mathrm{d}x}g\right) g^{-1} \right)^{\ast} ~ = ~ 
\left(g^{-1} \right)^{\ast}\left(\frac{\mathrm{d}}{\mathrm{d}x}g \right)^{\ast} ~ = ~ 
g\frac{\mathrm{d}}{\mathrm{d}x}g^{-1}
 ~ = ~ 
 -\frac{\mathrm{d}g}{\mathrm{d}x}g^{-1}
$$
$$
\left(gug^{-1}\right)^{\ast}
 ~ = ~ \left(g^{-1}\right)^{\ast}u^{\ast}g^{\ast}
  ~ = ~ -gug^{-1}
$$
and
$$
\left(ghg^{-1}\right)^{\ast}
 ~ = ~ \left(g^{-1}\right)^{\ast}h^{\ast}g^{\ast}
  ~ = ~ ghg^{-1}.
$$
Then (c) is easily derived from  
$$\omega= \frac{\mathrm{d}}{\mathrm{d}x}\left(a + a^\ast\right) - \left[a,~a^\ast\right] ~ = ~ 2\frac{\mathrm{d}}{\mathrm{d}x} h - \left[\left(h+u\right),~\left(h - u\right)\right] ~ = ~ 2\left(\frac{\mathrm{d}}{\mathrm{d}x} h - \left[u,~h\right]\right).$$
Equation (d) is the most interesting.

$$\delta h  ~ = ~ \frac{1}{2}  \delta\left(a  + a^\ast\right) ~ = ~ 
\frac{1}{2} \left(\left(\frac{\mathrm{d}}{\mathrm{d}x} \alpha + \frac{\mathrm{d}}{\mathrm{d}x} \alpha^\ast\right) - \left[a,\alpha\right] - \left[a,\alpha\right]^\ast\right) $$
$$ 
~ = ~
\frac{\mathrm{d}}{\mathrm{d}x} \alpha -\frac{1}{2}\left( \left[a,\alpha\right] - \left[a^\ast,\alpha^\ast\right]\right) ~ = ~ \frac{\mathrm{d}}{\mathrm{d}x} \alpha - \left[u,\alpha\right].$$

Likewise, 

 $$\delta u ~ = ~ \frac{1}{2}\left(\left(\frac{\mathrm{d}}{\mathrm{d}x} \alpha - \frac{\mathrm{d}}{\mathrm{d}x} \alpha^\ast\right) - \left[a,~\alpha\right] + \left[a,~\alpha\right]^\ast\right)$$
$$ 
~ = ~
- \frac{1}{2}\left(\left[a,~\alpha\right] + \left[a^\ast,~\alpha^\ast\right]\right) ~ = ~ -\left[h, ~\alpha\right].$$
\end{proof}

Now we set up the formalism for Yang-Mills on the torus.  We use the same notation for the connection one-forms.  The connection one-form is $\left(\frac{\mathrm{d}}{\mathrm{d}x} - u, ~\frac{\mathrm{d}}{\mathrm{d}y} +i h\right)$, the curvature is 
$$K ~=~ \left[\left(\frac{\mathrm{d}}{\mathrm{d}x} - u\right), ~\left(\frac{\mathrm{d}}{\mathrm{d}y} +i h\right)\right]
 ~=~ i\frac{\mathrm{d}h}{\mathrm{d}x}-\frac{\mathrm{d}}{\mathrm{d}y}u
-i\left[u,h\right],
$$  
and the (complex) gauge group  is just maps from the torus to $GL(n,\CC)$.  We define a theory on $\left[0,1\right]$ by assuming all the dependent variables do not depend on $y$.  We start on $\left[0,1\right]$, which we will double to $\left[-1,1\right]$.  We set up the theory on $\left[0,1\right]$ first.  Since we ultimately need some estimates, we also introduce function spaces. We use standard results from the Sobolev theory of connections and gauge transformations (see \cite{donaldson1997geometry},\cite{u82},\cite{frul12}).

\begin{prop}\label{prop3.1}
The dimensionally reduced Yang-Mills equations on $\left[0,1\right]$ can be given the following structure: 
\begin{itemize}
  \item[(a)] $\mA ~ = ~ \{(u, h):\left[0,1\right]\to \mathfrak{gl}(n,\CC) \mid \text{ both } u \text{ and } h \text{ are in the Hilbert space } H^1\}$.\\
  \item[(b)] $\tilde{\mG} ~ = ~ \{ g: \left[0,1\right]\to GL(n,\CC) \mid g(0) = g(1) = I,~ g\in H^2\}$.\\
  \item[(c)] $\tilde{\mU} ~ = ~ \{g: \left[0,1\right]\to U(n) \mid g\in H^2, g(0) = g(1) = I\}$.\\
  \item[(d)] The curvature $K ~ = ~ i \left[(\frac{\mathrm{d}}{\mathrm{d}x} - u),h\right]$ is smooth as a function in $H^0$ of  $\mA$ in $H^1$.\\  
  \item[(e)] The action of $\tilde{\mG}$ on $\mA$ is smooth.\\             
  \item[(f)] The action of $\tilde{\mU}$ on $\mA$ is an isometry.  Also the curvature $K$, up to the conventional, rather than our symplectic, identification of the Lie algebra with its dual, is the moment map for the action of $\tilde{\mU}$ on $\mA$.
\end{itemize}
\end{prop} 

Note that, as spaces of maps, $\tilde{\mG}=\mG^{\infty}$ and $\tilde{\mU}=\mU^{\infty}$. The different notation reflects their different origins in the two theories. Further, Proposition \ref{prop3.1} is just a restatement of what is in the literature.  The usual integral is
            $$||K||^2 ~ = ~ \left\langle K, K\right\rangle ~ = ~ \int_0^1  \Tr K K^\ast \mathrm{d}x.$$

There is no other natural function which is gauge invariant, but we can find two geometries on $\mA$ which give gradient flows reminiscent of the flows on $\mW_X$.  We could rework things from scratch explicitly on $\left[0,1\right]$.  However, after a slight adjustment of the gauge at the boundary, we can just double $\left[0,1\right]$ to $\left[-1, 1\right]$ with a symmetry and apply the theory for Yang-Mills on the torus.

We extend our tensors as follows.

\begin{df}\label{df3.2}
Extend $h$ by setting $h(-x) = h(x)$, $u(-x) = - u(x)$, $\left( X(-x)\right)^{\ast} = \left(X(x)\right)^{-1}$, and
$g(-x) = {g(x)^{-1}}^\ast$.  \\
Let 
$$\mA_\# = \{ \left(h, u \right)\from \left[-1,1\right] \to \gl\left(n,\CC\right) \text{ in } H^1 \mid h(x) = h(-x),~ u(-x) = - u(x) \}$$  
and let
$$\mG_\# = \{ g \from \left[-1,1\right] \to GL\left(n,\CC\right) \text{ in } H^2 \mid g(-x) = \left(g(x)^{-1}\right)^\ast \}.$$\\
\end{df}

\begin{prop}\label{prop3.4}
The theory for Yang-Mills on the torus $\frac{\left[-1,1\right]\times\left[-1,1\right]}{1=-1\times 1=-1}$ restricts to a consistent theory on the circle by having $a=u+h$ and $g$ be dependent only on $x$. 

In addition, if we set 
$$
a(x)= \left(a(-x)\right)^{\ast}$$
and
$$
g(x)^{\ast}= \left(g(-x)\right)^{-1}$$
we obtain a consistent theory on $\left[-1,1\right]$ (with $-1=1$). This is equivalent to a theory on $\left[0,1\right]$ with boundary conditions 
$$u(0)= u(1)= 0,$$
$$
\left(g(0)\right)^{\ast}= \left(g(0)\right)^{-1},$$
$$
\left(g(1)\right)^{\ast}= \left(g(1)\right)^{-1}.$$
\end{prop}

\begin{proof}
The reduction from $M\times\left[-1,1\right]$ to $M$ is standard. The extra connection $1$-form in the $y\in\left[-1,1\right]$ direction becomes the Higgs field, and the coupled theory is familiar.

The reduction from the circle $\frac{\left[-1,1\right]}{-1=1}$ to $\left[0,1\right]$ is not so familiar. The first step is to check that the conditions 
$g(x)^{\ast}= \left(g(-x)\right)^{-1}$ and $a(x)= a(-x)^{\ast}$ are consistent. But $g(x)^{\ast}= \left(g(-x)\right)^{-1}$ is consistent with the group action since
$$\left(g_1(x)g_2(x)\right)^{\ast}
~=~ g_2(x)^{\ast}g_1(x)^{\ast}
~=~ g_2(-x)^{-1}g_1(-x)^{-1}
 ~=~ \left(g_1(-x)g_2(-x)\right)^{-1}.
$$
The action of $g$ on $a$ is consistent since 
$$g(a)(x)^{\ast}
~=~ \left(\left(dg(x)+g(x)a(x)\right)g(x)^{-1}
\right)^{\ast}
~=~ \left(g(x)^{-1}\right)^{\ast}dg(x)^{\ast}+
\left(g(x)^{-1}\right)^{\ast}a(x)^{\ast}g(x)^{\ast}
$$
$$
~=~ -g(-x)d\left(g^{-1}\right)(-x)+g(-x)a(-x)g(-x)^{-1}
~=~ g(a)(-x).
$$

The only step in this chain which is confusing is 
$\left(dg(x)\right)^{\ast}= -d\left(g^{-1}\right)(-x)$ if
$g^{\ast}(x) = \left(g(-x)\right)^{-1}$.
The minus sign comes from the switch in derivative directions from $x$ to $-x$.

This theory clearly restricts to a theory on $\left[0,1\right]$ with the given conditions. However, if $a(0)^{\ast} = a(0)$ and $a(1)^{\ast} = a(1)$, noting that here $1=-1$, we can lift $a$ defined on $\left[0,1\right]$ to its double by setting $a(-x)=a(x)^{\ast}$. The same is true with $g(-x)=\left(\left(g(x)\right)^{-1}\right)^{\ast}$ if $g(0)=\left(\left(g(0)\right)^{-1}\right)^{\ast}$ and $g(1)=\left(\left(g(1)\right)^{-1}\right)^{\ast}$. So the restricted theory on  $\frac{\left[-1,1\right]}{-1=1}$ is equivalent to that theory on $\left[0,1\right]$.

\end{proof}

\begin{thm}\label{thm3.3}
 The two theories are equivalent via the unitary gauge groups:
 
\begin{itemize}
\item[(a)]	Every element in $\mA$ defined on $\left[0,1\right]$ is unitary gauge equivalent via an element in $\tilde{\mU}$ to a tensor in which $u(0) = u(1) = 0$ in $\mA$, which lifts to $\mA_\#$. \\
\item[(b)]	Every element in $\mA_\#$ restricts to an element of $\mA$.\\
\item[(c)]	Every element in $\tilde{\mG}$, after multiplication by an element in $\tilde{\mU}$, lifts to an element in $\mG_\#$.\\
\item[(d)] After the unitary adjustment in (c), the gauge group $\tilde{\mG}$ lifts to the subgroup $$\{ g(0) = g(1) = Id \}  \subset \mG_\#.$$
       The codimension of this subgroup in $\mG_\#$ is $2n^2$.\\
\item[(e)] The dimension of the moduli spaces is the same in each case.
\end{itemize}
\end{thm} 

\begin{proof}
Part (b) is easily seen to be true. An $H^1$ function on $\left[-1,1\right]$ restricts to an $H^1$ function on $\left[0,1\right]$, and elements in $\mA$ have no conditions at the boundary. 

For (a) note that $h(-x) = h(x)$ extends a function in $H^1$ to a function in $H^1$.  However, $u(-x) = -u(x)$ extends a function in $H^1$ to an $H^1$ function on the doubled circle only when $u(0) = u(1) = 0$. Now the (unitary) symmetry group acts on $u$ by $\left(g, u\right) \to \left(\frac{\mathrm{d}}{\mathrm{d}x} g + gu\right)g^{-1}$. One can get $u = 0$. However, such a $g$ would violate the boundary conditions.  On the other hand, it is possible to get $u(0) = u(1) = 0$ by choosing an essentially arbitrary smooth unitary $g$ with $g(0) = g(1) = I$, $\frac{\mathrm{d}}{\mathrm{d}x} g (0)= - u(0)$, and $\frac{\mathrm{d}}{\mathrm{d}x} g(1) = - u(1)$.  Note that each $u$ in $H^1$ is continuous, so $u(1)$ and $u(0)$ are well-defined. 

The condition that $u(-x)=-u(x)$ would mean that $u(0)=0$ and $u(1)=0$ are necessary for the reflection to be smooth. However, we could try another route. While there is no canonical way to do it, we can make a gauge transformation so that $u(0)=u(1)=0$. Then the reflection gives $u\in H^1$, as required. The flows are defined and the flows are invariant under gauge transformation, so we simply apply the gauge transformation, flow, and undo the gauge transformation to get to the result we want.

This points out something which might be important for numerical computations. If there is a discontinuity in the reflection, this discontinuity is preserved under the flows. With no discontinuity, one can take advantage of this regularity in setting up the numerics at the boundary. We are not experts in numerics, but think it worth thinking about.

For (c), let $g\in\tilde{\mG}$. The reflection $g(x)=(g(x)^{-1})^{\ast}$ is in $H^2$ if the derivatives agree at the reflection points. Since $g(0)=g(1)=I$, this requires $g_x(0)$ and $g_x(1)$ to be Hermitian. We can multiply $g$ by a smooth unitary $u\in\tilde{\mU}$ whose derivatives at 0 and 1 cancel the skew-Hermitian parts of $g_x(0)$ and $g_x(1)$. Then $ug$ has Hermitian derivatives at the endpoints and reflects to an element of $\mG_\#$. 

To see (d), note that if $g\in\mG_\#$, we have $g(0) = {g(0)^{-1}}^\ast$ and $g(1) = {g(1)^{-1}}^\ast$.  So $g$ is automatically in the unitary group at $\{0,1\}$. The dimension of the unitary group is $n^2$.  Hence the condition $g(0) = g(1) = I$ has codimension $2n^2$. 

Note that, by saying that after the unitary adjustment the lift of $\tilde{\mG}$ is smaller than $\mG_\#$ by $2n^2$ dimensions, 
(d) is also saying that the moduli spaces for our theory have dimension $2n^2$ larger than the moduli spaces for the restricted Yang-Mills theory.  However, the condition for lifting from $\mA$ to $\mA_\#$ also restricts to a submanifold of codimension $2n^2$. So the dimensions of the two theories coincide.

\end{proof}

The standard functional for Yang-Mills is the squared $L^2$-norm of the curvature. Under the identification above, the Yang-Mills curvature $K$ is proportional to the continuum moment map $\omega$. More precisely, if $q=\frac{\mathrm{d}}{\mathrm{d}x}h-\left[ u, h\right]$, then $K=iq$, while the moment map of Theorem \ref{prop2.1} is $\omega=2q$. Thus, their squared norms differ only by a constant factor, which only rescales the corresponding gradient flow.

  The inner product we defined on $\mA$ is the $L^2$ inner product, not the $H^1$ inner product:
         
$$\left\langle a_1,a_2\right\rangle ~ = ~ \int_0^1  \TrR  a_1 a_2^\ast~ \mathrm{d}x ~ = ~ \int_0^1\left( \TrR h_1 h_2 - \TrR u_1 u_2\right) \mathrm{d}x.$$\\
This means that the corresponding $L^2$-gradient does not define a smooth vector field on $\mA$, but its downward gradient flow is a diffusion (heat) equation. 
The Yang-Mills literature uses a different normalization of $a$, namely $a=2\left(h+u\right)$. This changes only the scaling of the flow parameter. We rewrite the equations using the normalization of this paper.
 
\begin{prop}\label{prop3.2}
The downward $L^2$ gradient flow on $\mA$ for $||K||^2=||q||^2$, with 
$q ~ = ~ \frac{\mathrm{d}}{\mathrm{d}x} h - \left[u, ~h\right]$, is

$$\frac{\mathrm{d}}{\mathrm{d}s}  h ~ = ~ 2~\left(\frac{\mathrm{d}}{\mathrm{d}x} q - \left[u,~ q\right]\right),$$

$$\frac{\mathrm{d}}{\mathrm{d}s} u  ~ = ~  -2~ \left[h, q\right].$$

The gradient flow on the circle and on elements satisfying the involution is gauge invariant and converges to 
$$M_0 ~ = ~ \left\{ \left(h, u\right)~ \,\middle|\, ~ \frac{\mathrm{d}}{\mathrm{d}x} h - \left[u, ~h\right] ~ = ~ 0\right\}.$$ 
\end{prop}

\begin{proof}   
By Theorem \ref{thm3.3} and the gauge invariance of the Yang-Mills flow, we can try to make the unitary gauge adjustment described above and apply the corresponding Yang-Mills theorem on a Riemann surface, applied to a torus with the symmetry we have described. The equation on the torus would be

$$\frac{1}{2} \frac{\mathrm{d}}{\mathrm{d}s} ih ~ = ~  ^\ast D_y \left(i q\right) ~ = ~  \frac{\mathrm{d}}{\mathrm{d}x} i q - \left[u, ~i q\right],$$

$$\frac{1}{2} \frac{\mathrm{d}}{\mathrm{d}s} u ~ = ~ ^\ast D_x \left(i q\right) ~ = ~ -\frac{\mathrm{d}}{\mathrm{d}y} i q + \left[i h, ~iq\right].$$

Since $q$ does not depend on $y$, this proves the formula for the flow. 
\end{proof} 

Unfortunately, the equations are actually a heat equation. Substituting $q=\frac{\mathrm{d}}{\mathrm{d}x}h-\left[u, ~h\right]$, we get
$$
\frac{1}{2} \frac{\mathrm{d}h}{\mathrm{d}s}
~ = ~ \left[\frac{\mathrm{d}}{\mathrm{d}x}-u, ~\left[\frac{\mathrm{d}}{\mathrm{d}x}-u, ~h\right]\right],
$$
$$
\frac{1}{2} \frac{\mathrm{d}u}{\mathrm{d}s}
~ = ~ -\left[h, ~\left[\frac{\mathrm{d}}{\mathrm{d}x}-u, ~h\right]\right].
$$

The problem is the heat flow is not defined for $a\in H^1$. One must go to $a\in H^2$. Then $\frac{\mathrm{d}}{\mathrm{d}x}a$ is continuous, and
$$
\frac{\mathrm{d}}{\mathrm{d}x}\left(a(-x)\right)
 ~ = ~ 
 -\frac{\mathrm{d}}{\mathrm{d}x}\left(a(-x)\right)
 ~ = ~ \left(\frac{\mathrm{d}}{\mathrm{d}x}\left(a(x)\right)\right)^{\ast}.
$$
Hence, 
$$-\frac{\mathrm{d}}{\mathrm{d}x}\left(a(-x)\right)
 ~ = ~ 
 \frac{\mathrm{d}}{\mathrm{d}x}\left(a(x)\right)\quad\text{or}$$ 
$$\frac{\mathrm{d}}{\mathrm{d}x}a(0) 
~ = ~ 
-\left(\frac{\mathrm{d}}{\mathrm{d}x}a(0)\right)^{\ast} 
.$$ 

We can assume $u(0)=0$ by gauge transformation. So the real extra condition is $\frac{\mathrm{d}}{\mathrm{d}x}\left(h(0)\right)=0$.
This mean $\omega(0)= \frac{\mathrm{d}}{\mathrm{d}x}\left(h(0)\right)
+\left[\omega(0), h(0)\right]=0$.
This is not a condition we want.

To do it by hand, we set
$$
\frac{\mathrm{d}}{\mathrm{d}s} \left\langle \omega,~ \omega \right\rangle 
 ~ = ~ 
 2\left\langle \omega,~ \frac{\mathrm{d} \omega}{\mathrm{d}s} \right\rangle 
  ~ = ~ 
4\left(
 \left\langle \omega,~ \left[\frac{\mathrm{d}}{\mathrm{d}x}-u,~
 \frac{\mathrm{d}h}{\mathrm{d}s}\right] \right\rangle
 -  \left\langle \omega,~ \left[\frac{\mathrm{d}u}{\mathrm{d}s},~h\right] \right\rangle
 \right)
 $$
 $$
  ~ = ~ 
 -4 \left\langle \left[\frac{\mathrm{d}}{\mathrm{d}x}-u,~ \omega\right],~
 \frac{\mathrm{d}h}{\mathrm{d}s} \right\rangle +
4 \left\langle  \left[h,~\omega\right],~  \frac{\mathrm{d}u}{\mathrm{d}s}
  \right\rangle,
  $$
 where we assume $\frac{\mathrm{d}h}{\mathrm{d}s}(0) =\frac{\mathrm{d}h}{\mathrm{d}s}(1)=0$ in order to integrate by parts.
 We then get the correct gradient
 $$
 \frac{\mathrm{d}h}{\mathrm{d}s}=
 -4\left[\left(\frac{\mathrm{d}}{\mathrm{d}x}-u\right),\omega\right],
 \quad\quad \quad\quad 
 \frac{\mathrm{d}u}{\mathrm{d}s} = 4\left[h, \omega\right],
 $$
 but with the extra condition $\frac{\mathrm{d}h}{\mathrm{d}s}(0)= \frac{\mathrm{d}h}{\mathrm{d}s}(1)= 0$.
 But then $\frac{\mathrm{d}}{\mathrm{d}s}\omega(0) = \frac{\mathrm{d}}{\mathrm{d}s}\omega(1) =0$ appears as an extra condition ($u(0)=u(1)=0$ can be assumed by changing gauges).

So the heat equation approach, the one that is used for Yang-Mills, is not viable. If we had started from scratch, we would have seen this right away.

\begin{prop}\label{propB'} The usual gradient descent using a heat equation does not restrict to our problem.
\end{prop}

Identifying a model element in $M_0$ for each orbit of $\tilde{\mU}$ turns out to be impossible. The best we can do is given in Theorem \ref{thm3.4}. We define the extended unitary gauge group as 
$$\tilde{\mU}^e ~ = ~ \left\{g\in H^2\left(\left[0,1\right], U(n)\right) ~ \,\middle|\, ~  g(0) = I\right\},$$ 
so that $\tilde{\mU} ~ = ~ \left\{g\in \tilde{\mU}^e ~ \,\middle|\, ~  g(1) = I\right\}$.

\begin{thm}\label{thm3.4}
 $M_0/\tilde{\mU}$ is of dimension $2 n^2$. Further, each $\left(h, u\right)\in M_0$ is gauge equivalent by $\tilde{\mU}^e$ to some $\left(h_0, 0 \right)$, where $h_0$ is a constant Hermitian matrix. If $g_j\in\tilde{\mU}^e$ take $\left(h_j, u_j\right)\in M_0$ to the same $\left(h_0, 0 \right)$, then they are $\tilde{\mU}$-gauge equivalent if and only if $g_1(1) = g_2(1)$. 
\end{thm}

\begin{proof} Here we assume $g$ is unitary.  For $q = \frac{\mathrm{d}}{\mathrm{d}x} h - \left[u,h\right]$, the real gauge group acts by  $u \mapsto \left(\frac{\mathrm{d}}{\mathrm{d}x} g + g u\right)g^{-1}$ and $h\mapsto g h g^{-1}$.  (This is compatible with $q \mapsto gq g^{-1}$, which is what it means that the theory is gauge invariant).

We simply solve $\frac{\mathrm{d}}{\mathrm{d}x} g + g u = 0$, $g(0) = I$. Since $u$ is skew-Hermitian, $g$ is unitary. Thus $g$ lies in $\tilde{\mU}^e$, with $g(1)$ arbitrary. In this case, $h\mapsto g h g^{-1}= h_0$ satisfies $\frac{\mathrm{d}}{\mathrm{d}x} g h g^{-1}= 0$. Thus, after dividing out by $\tilde{\mU}^e$, the remaining parameter is a constant Hermitian matrix $h_0 = g h g ^{-1}$ of dimension $n^2$.

If $g_j\in\tilde{\mU}^e$ take $\left(h_j, u_j\right)$ to the same  $\left(h_0, 0\right)$, then $g_2^{-1}g_1$ gives a gauge transformation from $\left(h_1, u_1\right)$ to $\left(h_2, u_2\right)$. The gauge transformation $g_2^{-1}g_1$ lies in $\tilde{\mU}$ if and only if $g_1(1) = g_2(1)$. Conversely, if $k\in\tilde{\mU}$ takes $\left(h_1, u_1\right)$ to $\left(h_2, u_2\right)$, then both $g_1$ and $g_2k$ are solutions to $\frac{\mathrm{d}}{\mathrm{d}x} g + g u_1 = 0$, $g(0) = I$. Thus $g_1 = g_2k$, and hence $g_1(1) = g_2(1)$. Therefore, an element of $M_0/\tilde{\mU}$ is determined by a constant Hermitian matrix $h_0$ and the endpoint $g(1)\in U(n)$. Each of these spaces has real dimension $n^2$. Thus $M_0/\tilde{\mU}$ has dimension $2n^2$.   
\end{proof}

\begin{cor}\label{cor3.2}
$\mA/\tilde{\mG} ~ = ~ M_0/\tilde{\mU}$.
\end{cor}

\begin{proof} 
The $\tilde{\mG}$-orbits are the fibers of the endpoint map $a\mapsto X_a(1)$. By the polar decomposition, each such orbit contains an element of $M_0$. Further, by Theorem \ref{thm3.4} and the uniqueness of polar decompositions, elements of $M_0$ in the same $\tilde{\mG}$-orbit are $\tilde{\mU}$-gauge equivalent. 
\end{proof}

\begin{cor}\label{cor3.4}
Every element in $M_0$ is $\tilde{\mU}$-gauge equivalent  to some $\left(\tilde{h},c\right)$ where $c$ is constant and $\tilde{h} = exp(cx) h_0 exp(-cx)$.
\end{cor}

\begin{proof} Choose $g_0\in\tilde{\mU}^e$ which takes $\left(h,u\right)$ to $\left(h_0, 0\right)$. Solve $g_0(1)^{-1} = e^c$ where $c$ is a constant. Then apply $g_c(x) = exp(cx)\in\tilde{\mU}^e$ to $\left(h_0,0\right)$. This gives an element in $M_0$ of the claimed form.  Then $g_c g_0$ gives the gauge transformation described in the corollary.
\end{proof} 

\begin{rmk}
To solve  $g_0(1) = exp(-c)$, note that $exp(-cx)$ is a geodesic in the unitary group from I to $exp(-c)$.  Since the unitary group is connected and compact, there is always a geodesic between any two points.  The difficulty is that there are many geodesics from I to $exp(-c)$ in the unitary group.  Even choosing the minimal one leaves one with ambiguities in the case that we are at a focal point. We suspect this must be true for finite $N$ as well.
\end{rmk}

\section{The geometry of the group orbits and regularizing flows}
	\label{s4:GeometryGrpOrbitsAndRegularizingFlows}  

We now come to looking at the regularizing flow. In order to understand things properly we need to look carefully at the situation for finite $N$. We have two geometries we can put on $\mW$, and two choices of regularizing functions, $\left\langle W, W\right\rangle$ and $\left\langle \Omega, \Omega\right\rangle$. We can identify each $\mW_Y$ with the gauge group $\mG$. There is no canonical identification, but the geometry is the same for every identification.  We will call this the group geometry.  We can identify $\mW_Y$ as a submanifold of $\mW$ which is an open set of a linear space. We can project any derivative at a point in $\mW_Y$ to be tangent to  $\mW_Y$. We call this the induced linear geometry.

The continuum flows introduced below are geometrically natural. We expect that they arise as limits of the corresponding projected flows in the finite-dimensional model, but do not yet know how to prove this (see Conjecture \ref{conj1}).

\begin{thm}\label{continuumfloweqs}  
The gradient flow for $\left\langle W, W\right\rangle$ in the group geometry and the gradient flow for $\left\langle \Omega, \Omega\right\rangle$ in the induced linear geometry are multiples of each other. Up to flow parameter rescaling, the flow is

          $$\frac{\mathrm{d}}{\mathrm{d}s} W_p  ~=~ -4 \left(-\Omega_p W_p + W_p \Omega_{p-1}\right).$$
\end{thm}

\begin{proof}[Proof of Theorem \ref{continuumfloweqs}]  
The gradient flow for $\left\langle \Omega,~\Omega\right\rangle$ in the induced linear geometry is the easiest to compute.

     $$\delta \left\langle \Omega,~\Omega\right\rangle
~ = ~ 
2 \sum_p  ~\left\langle ~\Omega_p,~  \left(W_{p+1}^{\ast} \delta W_{p+1} - \delta W_p W_p^{\ast}
+ \delta W_{p+1}^{\ast}W_{p+1} - W_p\delta W_p^{\ast}
\right)~\right\rangle
     $$
     $$
~ = ~ 
4 \sum_p  ~\left\langle ~\Omega_p,~  W_{p+1}^{\ast} \delta W_{p+1} - \delta W_p W_p^{\ast}~\right\rangle
.$$
 Continuing with the calculation yields:

$$\delta \left\langle \Omega,~\Omega\right\rangle~ = ~ 4 \sum_p \left\langle - \Omega_p W_p  + W_p \Omega_{p-1},~ \delta W_p\right\rangle.$$\\
In the induced linear geometry this means that

  $$\frac{\mathrm{d}}{\mathrm{d}s} W_p ~ = ~  4 \left(\Omega_p W_p - W_p \Omega_{p-1}\right).$$ \\
Note that this is automatically tangent to the orbit $\mW_X(1)$. We do not have to project.

The gradient flow for $\langle W,~W\rangle$ is slightly more complicated.  Here we write

    $$\delta \left\langle W,~W\right\rangle ~ = ~ 2 \sum_p \left\langle W_p,~ \delta W_p\right\rangle ~ = ~
2\sum_p \left\langle W_p,~\delta G_p G_p^{-1} W_p - W_p \delta G_{p-1}G_{p-1}^{-1}\right\rangle$$
$$~ = ~2 \sum_p \left\langle W_p W_p^{\ast} - W_{p+1}^{\ast}W_{p+1},~ \delta G_p G_p^{-1}\right\rangle
~ = ~ - 2 \sum_p \left\langle \Omega_p,~ \delta G_p G_p^{-1}\right\rangle.$$\\
In the group geometry, the gradient flow is

 $$\frac{\mathrm{d}}{\mathrm{d}s} G_p G_p^{-1} ~ = ~  2 \Omega_p.$$\\
But then we have

   $$\frac{\mathrm{d}}{\mathrm{d}s} W_p ~ = ~ \frac{\mathrm{d}}{\mathrm{d}s} G_p G_p^{-1} W_p - W_p \frac{\mathrm{d}}{\mathrm{d}s} G_{p-1} G_{p-1}^{-1}
~ = ~ 2 \left(\Omega_p W_p -  W_p \Omega_{p-1}\right).$$\\
These differ by a factor of 2.
\end{proof}

It seems reasonable to check what happens as $N \to \infty$.

\begin{df}\label{group_geometry_flow}
We call the flow
$$\frac{\mathrm{d}}{\mathrm{d}s} W_p
~ = ~ 2 \left(\Omega_p W_p -  W_p \Omega_{p-1}\right)$$ the \emph{group geometry flow}.
\end{df}
  
\begin{thm}\label{prop5.6}
The rescaled limit as $N \to \infty$ of the group geometry flow is
\begin{equation}\label{limitingeq}
\frac{\mathrm{d}}{\mathrm{d}s} a ~ = ~ \left[\left(\frac{\mathrm{d}}{\mathrm{d}x} - a\right),~\omega\right]\quad \text{ with }\quad \omega ~ = ~ \frac{\mathrm{d}}{\mathrm{d}x}\left(a + a^{\ast}\right) - \left[a,~a^{\ast}\right].
\end{equation}
\end{thm}

Before we go to the proof, it is appropriate to explain why this is a ``heat'' equation. It is an easy exercise to change variables.

\begin{cor}\label{cor5.7}
In terms of  $a = h + u$, we have $a^{\ast} = h - u$, and the gradient is given by

         $$\frac{\mathrm{d}}{\mathrm{d}s} h ~ = ~ 2~  \left[\frac{\mathrm{d}}{\mathrm{d}x} - u,~ \left[\frac{\mathrm{d}}{\mathrm{d}x} - u ,~ h\right]\right],$$
$$\frac{\mathrm{d}}{\mathrm{d}s} u ~ = ~  -2 ~ \left[h,~ \left[ \frac{\mathrm{d}}{\mathrm{d}x} - u,~h\right]\right].$$
\end{cor}

This is clearly a heat equation coupled to a lower order term, and the same flow as we got from Yang-Mills in Section \ref{s3:DimReducedYangMillsFlow}.

\begin{proof}[Proof of Theorem \ref{prop5.6}]  
Since the two flows in Theorem \ref{continuumfloweqs} differ only by an overall constant factor, we absorb this constant into the flow parameter and write
$$\frac{\mathrm{d}}{\mathrm{d}s} W_p ~ = ~ \Omega_p W_p - W_p \Omega_{p-1}.$$
We use the formulas we derived for the approximation term  $\left(1/N\right)^2\omega\left(p/N\right)$ which we found when we computed the limit $\Omega_p$ as $N \to \infty$.
    $$\frac{1}{N}  \frac{\mathrm{d}}{\mathrm{d}s} a\left(\frac{p}{N}\right)~ =~ 
    \frac{1}{N^3}\left( \frac{ \omega\left(\frac{p}{N}\right)-\omega\left(\frac{p-1}{N}\right)}{1~/~N}
+\omega\left(\frac{p}{N}\right)a\left(\frac{p}{N}\right) - a\left(\frac{p}{N}\right)\omega\left(\frac{p-1}{N}\right)\right).$$\\
Again, taking the plausible limit with $s^{\ast} = s /N^2$ gives

     $$\frac{\mathrm{d}}{\mathrm{d}s^{\ast}} a ~= ~ \frac{\mathrm{d}}{\mathrm{d}x}\omega - \left[a,~\omega\right].$$\\
Note that this agrees with the heat equation we got from the reduction of Yang-Mills. But for the same reason, we conjecture that it is not a viable flow.
\end{proof}

\begin{conjecture} The flows for the first finite dimensional flow
$$
\frac{\mathrm{d}}{\mathrm{d}s}W_p ~=~
2\left(\Omega_p W_p-W_p\Omega_{p-1}\right)
$$
do not converge to a viable flow as $N\to\infty$.
\end{conjecture}

The second flow is somewhat harder to describe, but is all we have left.  It is a gradient flow for $\left\langle W, W\right\rangle$ in the induced linear geometry, but it turns out to be also a gradient flow for $\left\langle\Omega, \Omega\right\rangle$ of an unusual sort. This is also proved in \cite{km25}, but it is key to our argument so we repeat the proof:

\begin{df}\label{linear_geometric_flow}
The \emph{linear geometric flow} is the flow induced by
$$\frac{\mathrm{d}}{\mathrm{d}s} W_p ~=~ -2\left( W_p - V_p\right),$$
where $V$ is the normal component of $W$ with respect to the group orbit. Namely,
$$W_p - V_p ~=~ \alpha_p W_p - W_p \alpha_{p-1},$$
$$ V_pW_p^{\ast} - W_{p+1}^{\ast} V_{p+1} ~=~ 0.$$ 

\end{df}
  
\begin{thm}\label{prop5.8}
The gradient flow in the geometry of the group orbit included in the linear space $\mathfrak{gl}(n,\CC)^N$ is the linear geometric flow.
\end{thm}
         
\begin{proof} 
Using the orthogonal decomposition of the ambient linear space into the tangent and normal spaces to $W_X$, we have for each $p=1,\dots, N$ that
     $$W_p - V_p ~=~ \alpha_p W_p - W_p \alpha_{p-1} $$
for some $\alpha_p$, where $V_p$ is normal to $W_X$. The fact that $V_p$ is normal means

        $$\sum_p \left\langle V_p, \left(\beta_p W_p - W_p \beta_{p-1}\right) \right\rangle~
      =~ \sum_p \left\langle \left(V_p W_p^{\ast} - W_{p+1}^{\ast} V_{p+1}\right), \beta_p\right\rangle~ =~ 0$$\\
for all $\beta=\left(\beta_1, \dots, \beta_{N-1}\right)$ where each $\beta_p\in\mathfrak{gl}(n,\CC)$ and $\beta_0=\beta_N=0$.  This is equivalent to $V_pW_p^{\ast} - W_{p+1}^{\ast}V_{p+1} = 0$.
\end{proof} 

To understand the impact of the flow on the imbalance in W, we consider the evolution of the moment map:

\begin{thm}\label{thm5.9}  
Along the gradient flow of Theorem \ref{prop5.8}, the finite-dimensional moment map $\Omega$ satisfies
$$\frac{\mathrm{d}}{\mathrm{d}s} \Omega ~=~ - 4 \Omega.$$
\end{thm} 

\begin{proof}    
Let 
$$H_p ~=~ W_{p+1}^{\ast}\frac{\mathrm{d}}{\mathrm{d}s}W_{p+1} - \frac{\mathrm{d}}{\mathrm{d}s}W_p W_p^{\ast}.$$\\
Then 

$$\frac{\mathrm{d}}{\mathrm{d}s} \Omega_p ~=~ H_p + H_p^{\ast}.$$

$$H_p ~=~2\left( W_{p+1}^{\ast}\left(-W_{p+1} + V_{p+1}\right) - \left(-W_p  + V_p\right)W_p^{\ast}\right)
     ~=~ -2 \Omega_p  +   2\left(W_{p+1}^{\ast}V_{p+1} - V_pW_p^{\ast}\right)  ~=~  -2 \Omega_p.$$\\
Since $\Omega_p^{\ast}=\Omega_{p}$, we also have $H_p^{\ast}=-2\Omega_{p}$. It follows that 
$$\frac{\mathrm{d}}{\mathrm{d}s} \Omega ~=~ -4 \Omega.$$
\end{proof} 

It was originally not clear how to take the continuum limit of the calculation in Theorem \ref{thm5.9}. However, the resulting equation suggests that we might instead look for a continuum flow satisfying 
$$\frac{\mathrm{d}}{\mathrm{d}s} \omega ~=~ -4 \omega.$$
This approach is somewhat similar to constructing a gradient descent flow for Yang-Mills, which is a vector field on an infinite dimensional manifold, not a heat equation. This is not the usual approach.

\begin{prop}\label{prop4.7new} If $\left[\frac{\mathrm{d}}{\mathrm{d}x}-a, \beta\right]=0$ and either $\beta(0)=0$ or $\beta(1)=0$, then $\beta=0$.
\end{prop} 

\begin{proof} 
Recall $\left(\frac{\mathrm{d}}{\mathrm{d}x}-a\right) X=0$. So $aX ~=~ \frac{\mathrm{d}}{\mathrm{d}x} X$ and $a ~=~ (\frac{\mathrm{d}}{\mathrm{d}x} X)X^{-1}$.
Hence,
$$\frac{\mathrm{d}}{\mathrm{d}x} \left(X^{-1} \beta X\right) ~=~  \left[ X^{-1}\beta X, ~X^{-1}\frac{\mathrm{d}}{\mathrm{d}x}X\right] + X^{-1}\left(\frac{\mathrm{d}}{\mathrm{d}x} \beta\right) X
$$
$$~=~ 
\left[ X^{-1}\beta X, ~X^{-1}aX\right] + X^{-1}\left(\frac{\mathrm{d}}{\mathrm{d}x} \beta\right) X 
~=~  X^{-1}\left(\left[\beta, a\right]+\frac{\mathrm{d}}{\mathrm{d}x}\beta\right)
 ~=~ 0.$$ 
Hence, $\frac{\mathrm{d}}{\mathrm{d}x} (X^{-1}\beta X) ~=~ 0$ and $X^{-1}\beta X$ is constant. Since $\beta(0)=0$ or $\beta(1)=0$, we have $X^{-1}(0)\beta(0) X(0)=0$ or $X^{-1}(1)\beta(1) X(1)=0$. Hence $X^{-1}\beta X$ being constant implies  $X^{-1}\beta X=0$, and $\beta=0$.
\end{proof} 

Proposition \ref{prop4.7new} is a standard fact. However, we put it in a proposition, for those unfamiliar with this theory.

\begin{lem}\label{lemConstructedFlow}  
Consider the tangent vector  
$$\delta h ~=~ \left[\left(\frac{\mathrm{d}}{\mathrm{d}x} - u\right),~\tau\right], \quad\quad\delta u = \left[\tau, h\right],$$ 
where $\tau ~=~ \tau^{\ast}$ and $\tau(0)=\tau(1)=0$. Then, recalling that $\omega ~=~ \left[\frac{\mathrm{d}}{\mathrm{d}x}-a, (a+a^{\ast})\right]$, we have that

    $$-2\omega ~=~  \left[\left(\frac{\mathrm{d}}{\mathrm{d}x} - u\right),~\left[\left(\frac{\mathrm{d}}{\mathrm{d}x} - u\right),~\tau\right]\right] + \left[h,\left[\tau,~h\right]\right] ~=~ L_{(h,u)}\left(\tau\right)$$\\
has a unique solution $\tau\in H^2$ for $\omega\in H^0$ and $a ~=~ h+u\in H^1$.  The map from $a \in H^1$ to $\delta a ~=~ \delta h + \delta u$ in $H^1$ is a smooth vector field on $\mA$. 

Moreover, under this flow,  $$\frac{\mathrm{d}}{\mathrm{d}s} \omega = -4 \omega.$$
\end{lem}

\begin{proof}    
It is a straightforward application of Proposition \ref{prop4.7new}, to see that $L_{(h,u)}$ can be smoothly inverted for $h,u\in H^1$. 
To see that it is negative with no kernel, use the fact that $\tau(0)=\tau(1)=0$ and integration by parts to get
$$\left\langle \tau,~L_{(h,u)}(\tau)\right\rangle ~=~ 
-\left\lVert\left[\left(\frac{\mathrm{d}}{\mathrm{d}x}-u\right), \tau\right]\right\rVert^2
-\left\lVert\left[\tau, h\right]\right\rVert^2.$$
If we set $L_{(h,u)}(\tau)=0$, we then get $\left[\left(\frac{\mathrm{d}}{\mathrm{d}x}-u\right),\tau\right] ~=~ 0$.
It follows that $\tau ~=~ 0$.
Standard elliptic theory then gives the unique $H^2$ solution and its smooth dependence on $h,u\in H^1$.  

Finally, we obtain $\frac{\mathrm{d}}{\mathrm{d}s} \omega =2L_{(h,u)}(\tau) = -4 \omega$ directly from the definitions.
\end{proof}

At first we were unable to show that the flow constructed in Lemma \ref{lemConstructedFlow} is the continuum limit as $N \to \infty$, and when we actually managed to construct the limit, it turned out to be different.

\begin{prop}\label{prop5.11}
If there is a limit as $N \to \infty$, it must satisfy
$$\frac{\mathrm{d}}{\mathrm{d}s} a ~=~ 2\left(v - a\right) ~=~ 2~ \left[\left(\frac{\mathrm{d}}{\mathrm{d}x} - a\right),~ \tau\right]
\quad
\text{where}~ \tau(0)=\tau(1)=0.
$$\\
Furthermore, $\left[\left(\frac{\mathrm{d}}{\mathrm{d}x} + a^{\ast}\right),~ \left(v + a^{\ast}\right)\right] ~=~ 0$.\\
\end{prop}

\begin{proof}    
As in our earlier computations, we start with $a\from \left[0,1\right]\to\gl\left(n,\CC\right)$, $v\from \left[0,1\right]\to\gl\left(n,\CC\right)$, and $\tau\from \left[0,1\right]\to\gl\left(n,\CC\right)$ satisfying $\tau(0)=\tau(1)=0$.

We set $W_p = I + \frac{1}{N} a\left(p/N\right)$, and $V_p = I + \frac{1}{N} v\left(p/N\right)$, and $\alpha_p = -\tau\left(p/N\right)$. The minus sign in front of $\tau\left(p/N\right)$ serves to identify the finite-dimensional gauge parameter $\alpha_p$ from Theorem \ref{prop5.8} with minus the continuum parameter $\tau$ of this proposition. The first identity makes sense, because if $W_p=I$ for all $p$, then $V_p=I$ for all $p$.

Substitute into $\frac{\mathrm{d}}{\mathrm{d}s}W_p= -2\left(W_p-V_p\right)$ to get
$$\frac{\mathrm{d}}{\mathrm{d}s}a=2\left(a_{p/N}-v_{p/N}\right).$$ 
Then 
$$W_p-V_p ~=~ -\left(\alpha_pW_p-W_p\alpha_{p-1}\right)$$ 
becomes
$$\left(I + \frac{1}{N} a\left(\frac{p}{N}\right)\right) - 
\left(I + \frac{1}{N}v\left(\frac{p}{N}\right)\right) ~=~
-\left(\tau\left(\frac{p}{N}\right)\left(I + \frac{1}{N} a\left(\frac{p}{N}\right)\right) -
\left(I + \frac{1}{N} a\left(\frac{p}{N}\right)\right)\tau\left(\frac{p-1}{N}\right)\right)$$

$$a\left(\frac{p}{N}\right)-v\left(\frac{p}{N}\right) ~=~
-\left(\frac{\tau\left(\frac{p}{N}\right) - \tau\left(\frac{p-1}{N}\right)}{1~/~N}+
\tau\left(\frac{p}{N}\right)a\left(\frac{p}{N}\right)-
a\left(\frac{p}{N}\right)\tau\left(\frac{p-1}{N}\right)\right).$$

In the limit we get 
$$a(x)-v(x) ~=~
-\frac{\mathrm{d}}{\mathrm{d}x}\tau\left(x\right)+
\left[\tau\left(x\right), a\left(x\right)\right] ~=~
-\left[\left(\frac{\mathrm{d}}{\mathrm{d}x}-a\right), \tau\right].
$$

Finally, $ V_pW_p^{\ast} - W_{p+1}^{\ast} V_{p+1} ~=~ 0$ becomes
$$
0 ~=~
\left(I + \frac{1}{N}v\left(\frac{p}{N}\right)\right)
\left(I + \frac{1}{N}a\left(\frac{p}{N}\right)\right)^{\ast}-
\left(I + \frac{1}{N}a\left(\frac{p+1}{N}\right)\right)^{\ast}
\left(I + \frac{1}{N}v\left(\frac{p+1}{N}\right)\right)
$$
$$
0 ~=~
v\left(\frac{p}{N}\right)a\left(\frac{p}{N}\right)^{\ast}-
a\left(\frac{p+1}{N}\right)^{\ast}v\left(\frac{p+1}{N}\right)
+ \frac{v\left(\frac{p}{N}\right)+a\left(\frac{p}{N}\right)
-a\left(\frac{p+1}{N}\right)^{\ast} - v\left(\frac{p+1}{N}\right)}{1~/~N}
$$

This in the limit is equivalent to 
$$
\left[v(x),a^{\ast}(x)\right] ~=~
\frac{\mathrm{d}}{\mathrm{d}x}v(x)+\frac{\mathrm{d}}{\mathrm{d}x}a^{\ast}(x),
$$
which is equivalent to what is in the proposition.
\end{proof}

\begin{conjecture}\label{conj1} The projected flow of Theorem \ref{prop5.8} converges, as $N \to\infty$, to the continuum flow described in Proposition \ref{prop5.11}.
\end{conjecture}

Of course, we now need to prove that the system of equations from Proposition \ref{prop5.11} has a solution when $a \in H^1$.

\begin{thm}\label{thm5.12}  
The system of equations

$$v - a ~=~  
 \left[\left(\frac{\mathrm{d}}{\mathrm{d}x} - a\right),~ \tau\right], $$\\
$$\left[\left(\frac{\mathrm{d}}{\mathrm{d}x} + a^{\ast}\right),~ \left(v + a^{\ast}\right)\right]  ~=~ 0,$$\\
with $\tau(0)=\tau(1)=0$, has a unique solution $\tau \in H^2$, $v \in H^1$ smoothly dependent on $a \in H^1$.
\end{thm}

\begin{proof}    
Let $b  ~=~ v + a^{\ast}$.   The equations now become 
 $$b - \left(a + a^{\ast}\right) ~=~
  \left[\left(\frac{\mathrm{d}}{\mathrm{d}x} - a\right), \tau\right], $$
$$
\left[\left(\frac{\mathrm{d}}{\mathrm{d}x} + a^{\ast}\right), b\right]
 ~=~ 0.
$$
These equations have the geometric meaning of $b$ being normal to the tangent plane at the connection $a$ to the orbit of the gauge group $\mG^{\infty}$, and
$-b + \left(a + a^{\ast}\right)$ being the projection onto this tangent space.

We will want $b ~=~ a + a^{\ast}+
  \left[\left(\frac{\mathrm{d}}{\mathrm{d}x} - a\right), \tau\right]$. So we will solve for $\tau$ first.
$$
0  ~=~ \left[\frac{\mathrm{d}}{\mathrm{d}x} + a^{\ast}, ~b\right]
  ~=~ \left[\frac{\mathrm{d}}{\mathrm{d}x} + a^{\ast}, ~(a + a^{\ast})\right]
  +  \left[\frac{\mathrm{d}}{\mathrm{d}x} + a^{\ast},~ \left[\frac{\mathrm{d}}{\mathrm{d}x} - a,~ \tau \right]~ \right].
$$

Recall $\omega ~=~ \left[\left(\frac{\mathrm{d}}{\mathrm{d}x} + a^{\ast}\right), \left(a + a^{\ast}\right)\right]$ is the moment map.
So $-\omega ~=~ \left[\frac{\mathrm{d}}{\mathrm{d}x} + a^{\ast}, \left[\left(\frac{\mathrm{d}}{\mathrm{d}x} - a\right), \tau\right]\right]$.

Consider the operator
$$
L_a\left(\tau\right) ~=~ 
\left[\left(\frac{\mathrm{d}}{\mathrm{d}x} + a^{\ast}\right),~  \left[ \left(\frac{\mathrm{d}}{\mathrm{d}x} - a\right),~ \tau\right]\right],
$$
related to that given in the proof of Lemma \ref{lemConstructedFlow}. $L_a\from H^2_0\to H^0$, where $H^2_0$ is the set of $\tau\in H^2$ with $\tau(0)=\tau(1)=0$. It is negative definite, since
$$\left\langle L_a(\tau),\tau \right\rangle ~=~
\left\langle \left[\left(\frac{\mathrm{d}}{\mathrm{d}x} + a^{\ast}\right), \left[\left(\frac{\mathrm{d}}{\mathrm{d}x} - a\right), \tau \right] \right],\tau \right\rangle ~=~
-\left\langle \left[\left(\frac{\mathrm{d}}{\mathrm{d}x} - a\right), \tau\right], \left[\left(\frac{\mathrm{d}}{\mathrm{d}x} - a\right), \tau \right] \right\rangle \leq 0.
$$

It has no kernel, since if $L_a(\tau)=0$,
$$
\left\langle \left[\left(\frac{\mathrm{d}}{\mathrm{d}x} - a\right), \tau\right], \left[\left(\frac{\mathrm{d}}{\mathrm{d}x} - a\right), \tau \right] \right\rangle ~=~ 0.
$$
By Proposition \ref{prop4.7new}, $\tau=0$.
We now simply set $b ~=~ a + a^{\ast}+
  \left[\left(\frac{\mathrm{d}}{\mathrm{d}x} - a\right), \tau\right]$.
  
  A slight modification of this argument shows that $v$ is differentiably dependent upon $a$.
\end{proof}

\begin{cor}\label{corofthm5.12}  
The equation $\frac{1}{2} \frac{\mathrm{d}}{\mathrm{d}s} a ~=~ v - a$ for $$
v - a ~=~ \left[\left(\frac{\mathrm{d}}{\mathrm{d}x} - a\right), \tau \right] \quad \quad \tau(0) ~=~ \tau(1) ~=~ 0
$$
$$
\left[\left(\frac{\mathrm{d}}{\mathrm{d}x} + a^{\ast}\right), \left(v+a^{\ast}\right) \right]=0
$$
is an ordinary differential equation in the group orbit of $\mG^{\infty}$ on  $\mA$. It has a smooth local solution.
\end{cor}

\begin{proof}    
This is an automatic consequence of Theorem \ref{thm5.12}. Ordinary differential equations behave the same way on infinite dimensional manifolds as they do on finite dimensional ones.  
\end{proof}

\begin{cor}\label{cor4.12}   
Let $\frac{1}{2} \frac{\mathrm{d}}{\mathrm{d}s} a ~=~ v - a$ as in the theorem. Then $\frac{\mathrm{d}}{\mathrm{d}s} \omega ~=~ -4~\omega$.
\end{cor}

\begin{proof}    
$$\delta\omega 
~=~ 
\delta \left(\frac{\mathrm{d}}{\mathrm{d}x}
\left(a+a^{\ast}\right) + \left[a^{\ast}, a\right]\right)
~=~ 
\frac{\mathrm{d}}{\mathrm{d}x}\left(\delta a + \delta a^{\ast}\right)
+ \left[\delta a^{\ast}, a\right] + \left[a^{\ast}, \delta a\right]
~=~ 
\left[\left(\frac{\mathrm{d}}{\mathrm{d}x}-a\right), \delta a^{\ast}\right] + \left[\frac{\mathrm{d}}{\mathrm{d}x} + a^{\ast}, \delta a\right].
$$
So
$$
\frac{\mathrm{d}}{\mathrm{d}s}\omega
~=~ 
\left[\left(\frac{\mathrm{d}}{\mathrm{d}x}+a^{\ast}\right), 
\left(2\left(v-a\right)\right)\right]
+ \left[\left(\frac{\mathrm{d}}{\mathrm{d}x}-a\right), 
2\left(v^{\ast}-a^{\ast} \right)\right]
$$
$$
~=~ 
-2\left( \left[\left(\frac{\mathrm{d}}{\mathrm{d}x}+a^{\ast}\right),
 \left(a+ a^{\ast}\right)\right]
 +  \left[\left(\frac{\mathrm{d}}{\mathrm{d}x}-a\right),
 \left(a+ a^{\ast}\right)\right]\right)
~=~ 
-4\omega.
$$

The proof is the same as for finite $N$.

$$\delta \omega ~=~  2 D + 2 D^{\ast}.$$  
Here
$$D ~=~ \frac{1}{2}\left[\left(\frac{\mathrm{d}}{\mathrm{d}x} + a^{\ast}\right), ~\delta a\right] ~=~ \left[\left(\frac{\mathrm{d}}{\mathrm{d}x} + a^{\ast}\right), ~\left(v-a\right)\right]$$
$$~=~ \left[\left(\frac{\mathrm{d}}{\mathrm{d}x} + a^{\ast}\right), ~\left(v + a^{\ast}\right)\right] ~-~  \left[\left(\frac{\mathrm{d}}{\mathrm{d}x} + a^{\ast}\right), ~\left(a + a^{\ast}\right)\right]
 ~=~  0 -  \omega.$$
\indent So 
$$\delta \omega ~=~ 2D + 2D^{\ast} ~=~ -4 \omega.$$
\end{proof} 

How should we interpret this result?  One would have thought that the natural flow  would be $\frac{\mathrm{d}}{\mathrm{d}s} a ~=~ -2\left(a + v\right)$ where   $\left(\frac{\mathrm{d}}{\mathrm{d}x} + a^{\ast}\right)v ~=~ 0$.

This projects the flow $\frac{\mathrm{d}}{\mathrm{d}s} a ~=~ - 2a$ onto the image of the operator $\tau\mapsto\left[\left(\frac{\mathrm{d}}{\mathrm{d}x} - a\right),~ \tau\right]$. However, we then get
$$\frac{\mathrm{d}}{\mathrm{d}s} \omega ~=~ -2~\left( \left[\left(\frac{\mathrm{d}}{\mathrm{d}x} - a\right), ~a^{\ast}\right] + ~\left[\left(\frac{\mathrm{d}}{\mathrm{d}x} + a^{\ast}\right),~ a\right]\right) ~=~ 
-2\omega + 2~\left[a,~ a^{\ast}\right].$$\\
\indent What we have done is to start with the flow $\frac{\mathrm{d}}{\mathrm{d}s} a  ~=~  -2\left(a + a^{\ast}\right)$ and
project it onto the image of $\tau\mapsto\left[\left(\frac{\mathrm{d}}{\mathrm{d}x} - a\right),~ \tau\right]$.  Then 
$$\frac{\mathrm{d}}{\mathrm{d}s} \omega ~=~ - 4 \omega.$$
It was not something we would have thought to do, but it in fact works.

We could easily show that the solution exists for all positive time due to the decay in $\omega$. However, what we need to know is how the coupled equation, which arises when we look at the gradients of the loss function plus the regularizing function, not just the regularizing function, behaves. This will use the same techniques that we use here, but we do that later.

\section{The interaction of the gradient flows}\label{s5}

The whole point of the construction is to find a minimum for $f\from GL(n,\CC)\to \RR$ by minimizing
          $$F(a) + \nu R(a)$$
where $F(a) = f(X(a))$ and $\nu R$ is a regularizing function which diverts the gradient flows toward the balanced locus $M_0 = \left\{a\mid \omega(a) = 0\right\}$.

We first compute the gradient flow for $F$.

\begin{prop}\label{prop5.1} The gradient flow for  $F$ (in the linear geometry) is given by

    $$\frac{\mathrm{d}}{\mathrm{d}s} a  ~ =~   -\left(X(x)^{-1}\right)^{\ast}~ X(1)^{\ast}~\left(\nabla_X f\right)~ X(x)^{\ast},$$
where $\nabla_X f$ is evaluated at $X=X(1)$.
\end{prop}

\begin{proof}
Using $\frac{\mathrm{d}}{\mathrm{d}x} X = aX$, it is an easy calculation, taking the derivative in $s$, to get

$$\delta \left( \frac{\mathrm{d}}{\mathrm{d}x} X X^{-1}\right) ~ =~ \delta a \quad\text{ and }  \quad    
  $$ 
$$\frac{\mathrm{d}}{\mathrm{d}x} \left(X^{-1} \delta X\right)
~ =~  
  - X^{-1}\frac{\mathrm{d}}{\mathrm{d}x} X^{-1} \delta X 
  + X^{-1}\frac{\mathrm{d}}{\mathrm{d}x} \left(\delta X\right)
  ~=~  
  X^{-1}  \delta \left(\left(\frac{\mathrm{d}}{\mathrm{d}x} X\right) X^{-1}\right)X
  ~ =~ 
  X^{-1} \delta a X.$$\\
So 

$$X(1)^{-1} \delta X(1) - X(0)^{-1} \delta X(0)~ =~ 
\int_0^1  X(x)^{-1}\delta a(x)X(x) \mathrm{d}x.$$\\
$\delta X(0) = 0$ because $X(0) = I$.

$$\delta F  ~ =~ \int_0^1 \left\langle\left(\nabla_X f\right),~ X(1) X(x)^{-1} \delta a X(x)\right\rangle \mathrm{d}x ~ =~
\int_0^1  \left\langle\left(X(x)^{-1}\right)^{\ast} X(1)^{\ast}\left(\nabla_X f\right) X(x)^{\ast},~ \delta a\right\rangle.$$

By the definition of the gradient flow (in the linear model) we have 
    $$\frac{\mathrm{d}}{\mathrm{d}s} a~ =~  -\left(X(x)^{-1}\right)^{\ast}~  X(1)^{\ast}~ \left(\nabla_X f\right)~  X(x)^{\ast}.$$
\end{proof}

In the applications of this formula, note that $Q = X(1)^{\ast}\left(\nabla_X f\right)$ does not depend on $x$.

\begin{lem}\label{lem5.2}  The bracket satisfies:
$$\left[\left(\frac{\mathrm{d}}{\mathrm{d}x} + a^{\ast}\right),~ \left(X(x)^{-1}\right)^{\ast} Q X(x)^{\ast}\right]~ =~ 0.$$
\end{lem}

\begin{proof} 
We know   
$$\frac{\mathrm{d}}{\mathrm{d}x}   \left(X^{(x)}\right)^{\ast}  ~ =~  \left(aX\right)^{\ast} ~ =~ X^{\ast}a^{\ast}$$ 
so 
$$\left(X^{\ast}\right)^{-1}\frac{\mathrm{d}}{\mathrm{d}x} X^{\ast} ~ =~ a^{\ast}.$$
Then
$$\left(\frac{\mathrm{d}}{\mathrm{d}x}\left(X^{\ast}\right)^{-1} Q X^{\ast}\right)  ~ =~ - \left[\left(X^{\ast}\right)^{-1}\left(\frac{\mathrm{d}}{\mathrm{d}x} X^{\ast}\right),~ \left(X^{\ast}\right)^{-1}Q X^{\ast}\right] ~ =~ -\left[ a^{\ast},~\left(X(x)^{-1}\right)^{\ast} Q X(x)^{\ast}\right].$$  
This cancels the term with $a^{\ast}$ and completes the proof.
\end{proof}

We note that the identity in Lemma \ref{lem5.2} is the limit of the corresponding finite-$N$ identities, although we don’t prove it here.

\begin{thm}\label{thm5.3}
The derivative of  $\left\langle\omega, \omega\right\rangle$ along the flow

$$\frac{\mathrm{d}}{\mathrm{d}s} a ~=~ -\left(X(x)^{-1}\right)^{\ast} ~X(1)^{\ast}~\left(\nabla_X f\right)~ X(x)^{\ast}$$ 
is 0.
\end{thm}

\begin{proof} 
  $\delta \left\langle\omega,\omega\right\rangle ~=~ 
  4 \left\langle \omega,~ \left[\left(\frac{\mathrm{d}}{\mathrm{d}x} + a^{\ast}\right), \delta a\right]\right\rangle$.

Because $\omega = \omega^{\ast}$, the variation with the terms $\delta a^{\ast}$ is equal to the variation of the terms with $\delta a$. This is because $\left\langle\omega, \beta\right\rangle ~=~ 
\left\langle\omega^{\ast}, \beta^{\ast}\right\rangle ~=~ 
\left\langle\omega, \beta^{\ast}\right\rangle$ for any $\beta$. Hence the 4 and the lack of  terms with $\delta a^{\ast}$. We plug in the formula for $\delta a = \frac{\mathrm{d}}{\mathrm{d}s} a$, which gives us a chance to use Lemma \ref{lem5.2}:
$$\frac{\mathrm{d}}{\mathrm{d}x} \frac{\mathrm{d}}{\mathrm{d}s} a - \left[\frac{\mathrm{d}}{\mathrm{d}s}a,a^{\ast}\right] ~= ~ \left[\left(\frac{\mathrm{d}}{\mathrm{d}x} + a^{\ast}\right), ~\frac{\mathrm{d}}{\mathrm{d}s} a\right] ~=~
-\left[\left(\frac{\mathrm{d}}{\mathrm{d}x} + a^{\ast}\right), (X(x)^{-1})^{\ast} Q X(x)^{\ast}\right] ~=~ 0.$$
\end{proof}

We note that this theorem is true for finite $N$. However, when using $\left\langle W, W\right\rangle$ as the regularizing function for finite $N$, we were unable to show the equivalent theorem is true. There seem to be a lot of reasons to use $\left\langle \Omega,\Omega\right\rangle$ as a regularizing functional for finite $N$.

Theorem \ref{thm5.3} shows the gradient flow for $F$ leaves $\left\langle \omega, \omega\right\rangle$ fixed. We now turn to the effect of the regularizing flows for $\left\langle \omega, \omega\right\rangle$ on $F$.  These flows are all tangent to the gauge orbits $\mA_X(1)$.  We use the parameter $t$, as these flows are distinct from the flow $\frac{\mathrm{d}}{\mathrm{d}s} a$ we used above.

\begin{thm}\label{thm5.4}
Suppose that the flow $\frac{\mathrm{d}}{\mathrm{d}s} a$ is tangent to  $\mA_X(1)$. Then $F$ is constant along the flow.
\end{thm}

\begin{proof}  
Along the flow,
$$ \frac{\mathrm{d}}{\mathrm{d}t} F ~=~  \left\langle \nabla_X(1) f ,~\frac{\mathrm{d}}{\mathrm{d}t} X(1)\right\rangle.$$  
But $\frac{\mathrm{d}}{\mathrm{d}t} X(1) = 0$ if the flow is tangent to $\mA_X(1)$.
\end{proof}

We are now ready to prove that the regularizing flows for $\left\langle \omega, \omega\right\rangle$  have the correct behavior.

\begin{thm}\label{thm5.5new} 
Suppose
$$
\frac{\mathrm{d}}{\mathrm{d}s} a ~=~ 
-\left(X(x)^{-1}\right)^{\ast}X(1)^{\ast}
\left(\nabla_Xf\right)X(x)^{\ast}+
\nu\frac{\mathrm{d}}{\mathrm{d}s} a,
$$
where the first term is the downward gradient for $F(a)$ and the second term is any downward gradient flow $\left\langle \omega, \omega \right\rangle$ tangent to the submanifold $X(1)$ fixed. Then
$$
\frac{\mathrm{d}}{\mathrm{d}s}F(a) ~=~ 
-\left\langle dF, \nabla F \right\rangle
 ~=~ 
-\left\langle \nabla F, \nabla F \right\rangle,
$$
$$
\frac{\mathrm{d}}{\mathrm{d}s}\left\langle \omega, \omega \right\rangle ~=~ 
\nu~ 2\left\langle \omega, \frac{\mathrm{d}\omega}{\mathrm{d}s} \right\rangle.
$$
\end{thm}

\begin{proof} 
Since we require $\frac{\mathrm{d}}{\mathrm{d}s}a$ to be tangent to the manifold $X(1)=constant$, the first term is automatic. The second term follows from Theorem \ref{thm5.3}.
\end{proof} 

\begin{cor}  
Let $\frac{\mathrm{d}}{\mathrm{d}s}a=2\left(v-a\right)$, where $v\in H^1$ is a smooth function of $a\in H^1$ satisfying 
$$v-a ~ = ~  
\left[ \frac{\mathrm{d}}{\mathrm{d}x}-a, \tau \right],\quad
\tau(0)=\tau(1)=0\quad \text{with}\quad
\tau\in H^2,\quad \text{ and}\quad$$
$$\left[\frac{\mathrm{d}}{\mathrm{d}x}+a^{\ast}, 
\left(v+a^{\ast}\right)
\right] ~ = ~ 0.$$
Then, under the hypotheses of Theorem \ref{thm5.5new},
$$
\frac{\mathrm{d}}{\mathrm{d}s} F ~ = ~ 
-\left\langle\nabla F,~ \nabla F\right\rangle,
$$
$$
\frac{\mathrm{d}}{\mathrm{d}s} \omega^2 ~ = ~ 
2\nu \left\langle\omega,~ \frac{\mathrm{d}}{\mathrm{d}s}\omega\right\rangle ~ = ~ 
-8\nu \left\langle\omega,~ \omega\right\rangle.
$$
\end{cor} 

\begin{proof}  
This is a straightforward application of Theorem \ref{thm5.5new} and Corollary \ref{cor4.12}.
\end{proof}

\bibliographystyle{alpha}

\bibliography{refs}

\end{document}